\documentclass{agtart_a}
\pdfoutput=1
\usepackage{pinlabel}

\title{Widths of surface knots}

\author{Yasushi Takeda}
\givenname{Yasushi}
\surname{Takeda}
\address{Graduate School of Mathematics\\
Kyushu University\\\newline
Hakozaki\\Fukuoka 812-8581\\Japan}
\email{takeda@math.kyushu-u.ac.jp}
\urladdr{}

\volumenumber{6}
\issuenumber{}
\publicationyear{2006}
\papernumber{64}
\startpage{1831}
\endpage{1861}

\doi{}
\MR{}
\Zbl{}


\received{7 February 2006}
\revised{10 August 2006}
\accepted{21 August 2006}
\published{1 November 2006}
\publishedonline{1 November 2006}
\proposed{}
\seconded{}
\corresponding{}
\editor{CPR}
\version{}

\arxivreference{}

\let\xysavmatrix\xymatrix
\def\xymatrix{\disablesubscriptcorrection\xysavmatrix}
\AtBeginDocument{\let\bar\wbar\let\tilde\wtilde\let\hat\what
\renewcommand{\setminus}{\smallsetminus}}
\def\psfrag#1#2{\relax}

\newcommand{\RP}{\mathbb{R}$P$}

\makeatletter
\def\cnewtheorem#1[#2]#3{\newtheorem{#1}{#3}[section]
\expandafter\let\csname c@#1\endcsname\c@theorem}
\makeatother

\newtheorem{theorem}{Theorem}[section]
\cnewtheorem{lemma}[theorem]{Lemma}
\cnewtheorem{proposition}[theorem]{Proposition}
\cnewtheorem{corollary}[theorem]{Corollary}
\cnewtheorem{conjecture}[theorem]{Conjecture}

\theoremstyle{definition}
\cnewtheorem{definition}[theorem]{Definition}
\cnewtheorem{example}[theorem]{Example}
\cnewtheorem{xca}[theorem]{Exercise}

\theoremstyle{remark}
\cnewtheorem{remark}[theorem]{Remark}

\numberwithin{equation}{section}

\begin{document}

\begin{abstract}
We study surface knots in 4--space by using generic planar projections.
These projections have fold points and cusps as their singularities 
and the image of the singular point set divides the plane into several
regions.
The width (or the total width) of a surface knot is a numerical invariant
related to the number of points in the inverse image of a point in
each of the regions. 
We determine the widths of certain surface knots and characterize
those surface knots with small total widths.
Relation to the surface braid index is also studied. 
\end{abstract}

\maketitle

\section{Introduction}\label{sec-1}

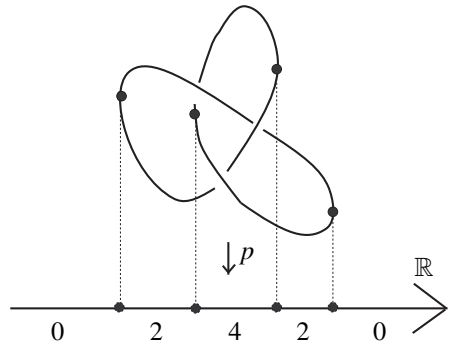
\begin{figure}[ht!]
\labellist\hair5pt\small
\pinlabel $\R$ at 442 218
\pinlabel {$0$} [t] at -89 156
\pinlabel {$0$} [t] at 378 156
\pinlabel {$2$} [t] at 54 156
\pinlabel {$4$} [t] at 168 156
\pinlabel {$2$} [t] at 267 156
\pinlabel {$p$} [l] at 158 234
\endlabellist
\begin{center}
\begin{minipage}{6cm}
\includegraphics[width=6cm]{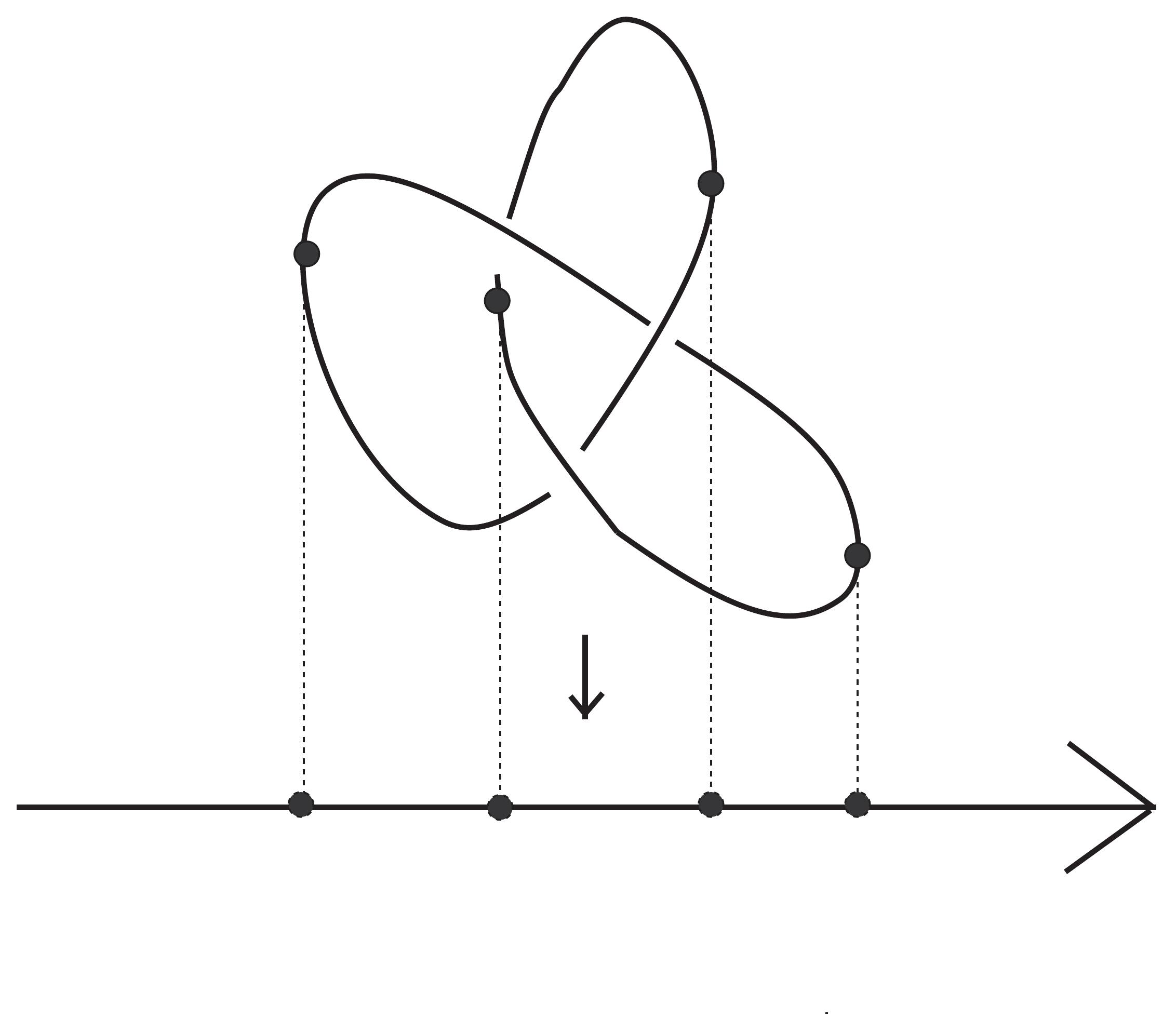}
\end{minipage}
\end{center}
\caption{Local widths of an embedded circle in $\R^3$} 
\label{fig:a}
\end{figure}

The notion of {\it width\/} for classical knots was introduced by
Gabai \cite{G} as a generalization of the bridge index, which
plays an important role in the classical knot theory. The width
was useful for solving difficult problems. 
More precisely, we consider a generic projection $p$ of 
an embedded circle in $\R^3$ into the line $\R$ as in
\fullref{fig:a}.
Then non-degenerate critical points appear as its singularities
and their images divide the line into several intervals.
For each such interval, we consider the number of points in 
$p^{-1}(x)$ for a point $x$ in the interval, and we call it
the \textit{local width}, which does not depend on the choice
of $x$.
The width of a knot is the minimum of the total of local widths 
over all embedded circles representing the given knot.

By a {\it surface knot\/},
we mean (the isotopy class of) a closed connected (possibly non-orientable)
smoothly embedded surface
in $\R^4$. For a surface knot, Carter--Saito \cite[Section 4.6]{cs2} considered
the analogy of the width.
They applied the notion of chart for 
the definition of width for surface knots.  
A chart is a planar projection of a surface knot together with
an associated graph, which was first introduced in the surface braid theory
(see Kamada \cite{Ka}).  
The graph is constructed by using a generic 
projection into $3$--space of a surface knot.  
The generic projections into $3$--space of surface knots
have double points, triple points, and branch points
as their singularities, and the 
charts represent the state of the combination
of these singularities.
Moreover, charts form several planar regions which are surrounded by curves
representing double points and fold lines, and the width of a surface
knot which Carter--Saito introduced is defined by using the number of
points in the fiber over a point (local width) in each of these
regions like the width for classical knots.
They considered the minimum (over all representatives of the given isotopy 
class) of the maximum of local widths over all the regions.

However, the width which Carter--Saito defined is slightly different from
the one which Gabai defined.
In fact, Carter--Saito considered the maximum of local widths for the
definition of width and Gabai considered the total of local widths.
Moreover, the width of surface knots has not been studied so much until now
as far as the author knows.

In this paper, for surface knots, we study the width defined by
Carter--Saito, and the {\it total width\/} which is the straightforward
analogy of the width for classical knots defined by Gabai.  For this
purpose, we consider generic planar projections of surface knots
instead of charts.  In the surface knot theory, we often use generic
projections into $3$--space: in fact, many results have been obtainted
by using projections into $3$--space, and since we can view the
diagrams in $3$--space, they facilitate the study of surface knots.
Generic planar projections have also been useful (for example, see
Carrara, Carter and Saito \cite{ccs}, Carrara, Ruas, and Saeki
\cite{crs}, Saeki and Takeda \cite{st} and Yamamoto \cite{Ya}).  Planar
projections have fold points and cusps as their singularities.  Cusps
appear as discrete points and fold points appear as a $1$--dimensional
submanifold.  Let us call the set of cusps and fold points in the
surface the {\it singular set\/}.  For a given surface knot, the image
of the singular set divides the plane into several regions.  For each
such region, we consider the number of points in the pre-image of a
point in that region and the maximum or the total of these numbers
over all the regions.  Then we take the minimum of these numbers over
all embedded surfaces representing the given surface knot.  Roughly
speaking, this defines the width and the total width of a surface
knot.

The paper is organized as follows. 
In \fullref{sec-2} we define the width and the total width of surface knots
and recall the definitions of the genericity of mappings and the
triviality of surface knots.
In \fullref{sec-3} we study the width and determine the width of some surface
knots such as ribbon surface knots and $n$--twist spun $2$--bridge knots.
In \fullref{sec-4} we consider the relationship between the width 
and the surface braid index and show that the width is always smaller
than or equal to the twice of the surface braid index plus two. 
We also show that in general the difference between these two invariants can be
arbitrarily large. 
In \fullref{sec-5} we give some characterization theorems of surface knots
with small total widths.

Throughout the paper,
we work in the smooth category. 

The author would like to thank Professor Osamu Saeki for
helpful suggestions and Professor Mitsuyoshi Kato for his 
constant encouragement.
He also thanks the referee for careful reading and useful comments.
The author has been supported by JSPS Research  
Fellowships for Young Scientists.
\section{Preliminaries}\label{sec-2}

In this section, we prepare several notions from singularity theory
and define the width and the total width for surface knots in $\R^4$.
For singularity theory, the reader is referred to Golubitsky and
Guillemin \cite{GG}, for example.

\begin{definition}
Let $F$ be a closed connected surface. 
Denote by $C^{\infty}(F, \R^2)$ the set of all smooth mappings from
$F$ to $\R^2$, endowed with the Whitney $C^{\infty}$ topology.
Let $f$ and $g$ be elements of $C^{\infty}(F, \R^2)$.
Then $f$ is \textit{equivalent} to $g$ if there exist diffeomorphisms
$p \co F \to F$ and $q \co \R^2 \to \R^2$ such that $q \circ f = g \circ p$.
\end{definition}

\begin{definition}
Let $f$ be an element of $C^{\infty}(F, \R^2)$.
Then $f$ is said to be \textit{$C^{\infty}$ stable} if there exists a
neighborhood $N_{f}$ of $f$ in $C^{\infty}(F, \R^2)$
such that each $g$ in $N_{f}$ is equivalent to $f$.
\end{definition}

\begin{definition}
Let $f \co F \to \R^2$ be a smooth mapping from $F$ to $\R^2$.
Then $q \in F$ is called a \textit{fold point} if
we can choose local coordinates $(x,y)$ centered at $q$
and $(U,V)$ centered at $f(q)$ such that $f$, in a neighborhood of $q$,
is of the form:
$$U=x,\ V=y^2.$$
Moreover, $q \in F$ is called a \textit{cusp} if we can choose
local coordinates $(x,y)$ centered at $q$ and $(U,V)$ centered at $f(q)$
such that $f$, in a neighborhood of $q$, is of the form:
$$U=x,\ V=xy+y^3.$$
\end{definition}

We denote by $S_{1}(f)$ the set of fold points and cusps,
and by $S^2_{1}(f)$ the set of cusps.

Note that $S_{1}(f)$ is a regular $1$--dimensional submanifold of $F$ and
$S^2_{1}(f)$ is a discrete set.

Recall the following well-known characterization of $C^{\infty}$
stable mappings in\break $C^{\infty}(F, \R^2)$.

\begin{proposition}
Let $f \co F \to \R^2$ be a smooth mapping from a closed connected
surface $F$ to $\R^2$.
Then $f$ is $C^{\infty}$ stable if and only if $f$ has only
fold points and cusps as its singularities, its restriction to the
set of fold points is an immersion with normal crossings, and
for each cusp $q$, we have:
 $$f^{-1}(f(q)) \cap S_{1}(f)=\{q\}.$$
\end{proposition}

Let $F$ be a closed connected surface. For a smooth map 
$f \co F \to \R^2$, we set
$$S(f) = \{x \in F \,|\ {\rm rank}\ df_{x} < 2\},$$
which is called the {\it singular point set} of $f$.
If $f$ is $C^{\infty}$ stable, then we clearly have
$S_{1}(f) = S(f)$.

The following theorem is well-known (see Thom \cite{Thom}).

\begin{theorem}
Let $f \co F \to \R^2$ be a $C^{\infty}$ stable mapping from a
closed connected surface $F$ to $\R^2$. 
Then the number of cusps of $f$ has the same parity as the Euler
characteristic $\chi (F)$ of $F$.

\end{theorem}

\begin{definition}
Let $f \co F \to \R^4$ be an embedding of a closed connected
surface.
Let $\pi \co \R^4 \to \R^2$ be an orthogonal projection.
Then we say that $\pi$ is \textit{generic} with respect
to $f$ (or with respect to $f(F)$) if $\pi \circ f$ is
$C^{\infty}$ stable.
\end{definition}

By Mather \cite{Ma}, almost every orthogonal projection is generic
with respect to $f$.

\begin{definition} 
Let $f \co F \to \R^4$ be an embedding of a closed connected surface
$F$.
Let $\pi \co \R^4 \to \R^2$ be an orthogonal projection 
which is generic with respect to $f$. 
In this cace,
$\pi \circ f$ has fold points and cusps as its singularities.
Let $S(\pi \circ f)$($\subset F$) denote the set of these singularities.
The singular value set $\pi \circ f(S(\pi \circ f))$ divides the plane
$\R^2$ into several regions.
For a point $x$ in a given region, we call the number of elements in
the set $(\pi \circ f)^{-1}(x)$ 
the {\it local width\/}, which does not depend on the choice of $x$
and is always even (see the proof of \fullref{lem:t1}).
Let $w(f,\pi)$ (or $w(f(F),\pi)$) be the maximum of the local widths 
over all the regions
and $tw(f,\pi)$ (or $tw(f(F),\pi)$) be the total of the local widths 
over all the regions.
The {\it width\/} $w(f(F))$ of a surface knot $f(F)$ is
the minimum of $w(\tilde{f},\tilde{\pi})$, where $\tilde{f}$ runs over
all embeddings isotopic to $f$ and $\tilde{\pi}$ runs over all orthogonal
projections which are generic with respect to $\tilde{f}$.
Moreover, the {\it total width\/} $tw(f(F))$ of a surface knot $f(F)$
is the minimum of $tw(\tilde{f},\tilde{\pi})$,
where $\tilde{f}$ runs over all embeddings isotopic to $f$ and 
$\tilde{\pi}$ runs over all orthogonal projections which are generic
with respect to $\tilde{f}$.
\end{definition}

Let us now recall  the definitions of a handlebody, the standard 
projective planes in $\R^4$ and  the normal Euler number. 

An orientable \textit{handlebody} is a compact orientable $3$--manifold
obtained by attaching a finite number of $1$--handles to a $3$--ball
(the number of $1$--handles may possibly be zero).
A non-orientable \textit{handlebody} is a compact non-orientable
$3$--manifold obtained by attaching a finite number of $1$--handles
to a $3$--ball.

\begin{figure}[ht!]
\labellist\hair5pt\small
\pinlabel {$\R^3 \times \{0\}$} [tr] at -84 315
\pinlabel {$D^2 \subset \R^3 \times [0,\infty)$} [bl] at 196 904
\pinlabel {attach} [l] at 268 609
\endlabellist
\begin{center}
\begin{minipage}{8cm}
\includegraphics[width=8cm]{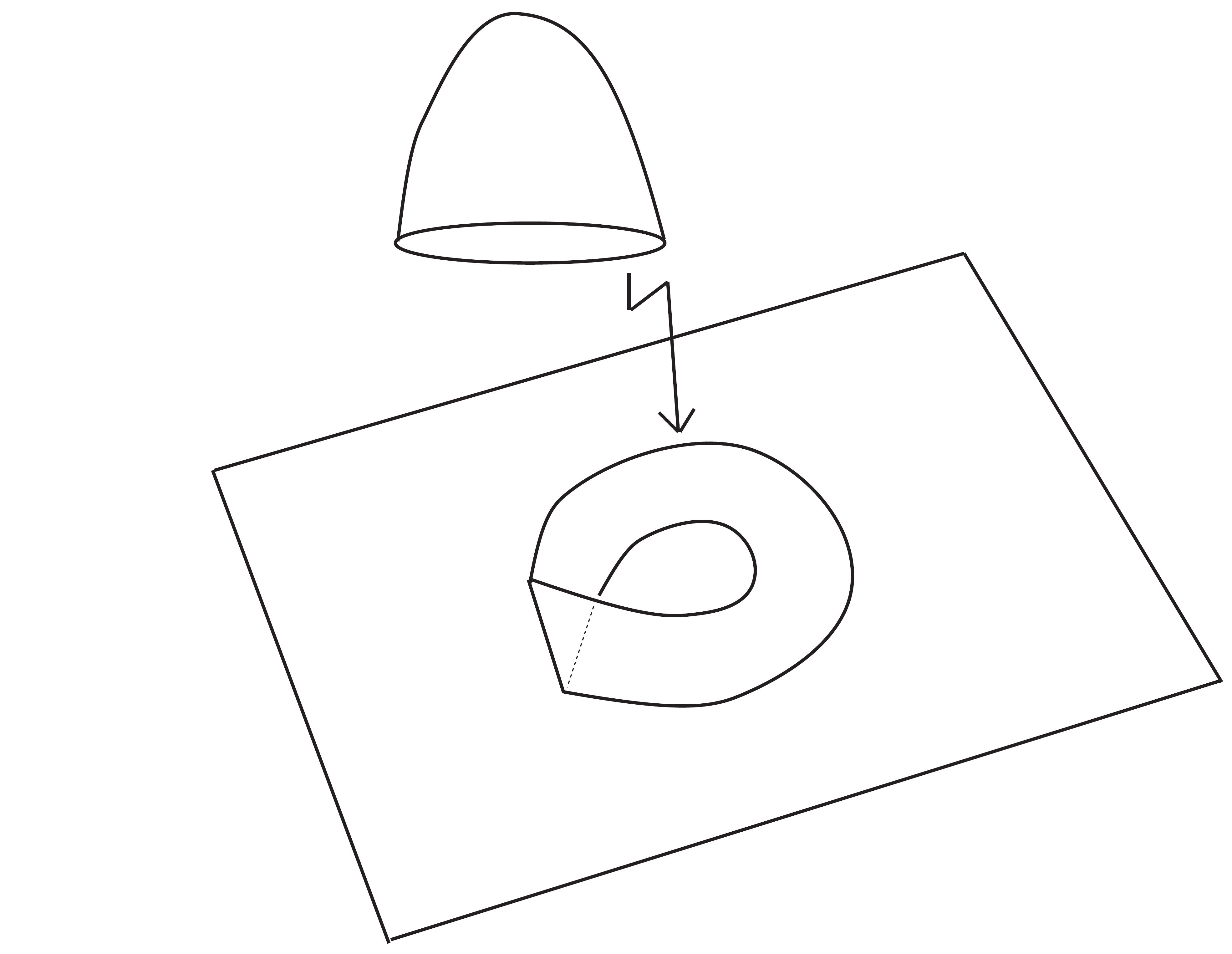}
\end{minipage}
\end{center}
\caption{The standardly embedded projective plane in $\R^4$}
\label{fig:b}
\end{figure}

The standardly embedded projective plane in $\R^4$
is constructed as in \fullref{fig:b}, by
attaching an unknotted disk in $\R^3 \times [0, \infty)$
to a ``trivially embedded" ${\rm M\ddot{o}bius}$ band in $\R^3 \times \{0\}$.
We have two trivially embedded ${\rm M\ddot{o}bius}$ bands up to isotopy, and
accordingly we have two kinds of standard projective planes in $\R^4$.
These surface knots have normal Euler number $\pm 2$.
Normal Euler number is an isotopy invariant of surface knots
(for example, see Carter and Saito \cite{cs2}).    

There are several definitions of trivial surface knots in the
litterature (for example, see Hosokawa and Kawauchi \cite{Ho}). 
In this paper, we adopt the following definition.

\begin{definition}
For a surface knot, we say that it is {\it strongly trivial\/} if
it is the boundary of a handlebody embedded in $\R^4$.
Moreover, we say that a non-orientable surface knot
is {\it trivial\/} if it is the connected 
sum of some copies of the standardly embedded projective planes
in $\R^4$, that is, the connected sum of $k$ copies of the standardly
embedded projective plane with normal Euler number $+2$ 
and $l$ copies with normal Euler number $-2$ for some $k \geq 0$ and
$l \geq 0$ with $k+l \geq 1$. 
\end{definition}

A surface knot is trivial if it is strongly trivial.
However, a trivial surface knot may not necessarily be strongly
trivial. In fact, if a surface knot is strongly trivial,
then its Euler characteristic must be even. 
More precisely, a trivial surface knot is strongly trivial
if and only if its normal Euler number vanishes.
Furthermore, for a closed connected non-orientable surface of
non-orientable genus $g$, the number of trivial surface knots
diffeomorphic to it is equal to $g+1$, and if $g$ is even,
then a strongly trivial surface knot diffeomorphic to it
exists and is unique (for example, see \cite{Ho}).  

The following lemma is often used throughout this paper.  

\begin{lemma}[Carrara, Ruas and Saeki~\cite{crs}]\label{lem:crs}
Let $\pi \co \R^4 \to \R^2$ and $\pi_1 \co \R^2 \to \R$
be orthogonal projections. 
For an embedding $f \co F \to \R^4$ of a closed connected surface
$F$,
if $\pi \circ f$ is $C^{\infty}$ stable without
cusps and $\pi_1 \circ \pi \circ
f \co F \to \R$ is a Morse function\footnote{A smooth function
on a smooth manifold is a \textit{Morse function} if its critical
points are all non-degenerate.}
 with at most four critical
points, then $f(F)$ is strongly trivial. 
\end{lemma}

\section{Widths of certain surface knots}\label{sec-3}

In this section, we characterize those surface knots with width two
and determine the widths of ribbon surface knots and 
$n$--twist spun $2$--bridge knots. 

Let we begin by the following lemma.

\begin{figure}[ht!]
\labellist\small
\pinlabel {$D^2_{i}$} [b] at -6 544
\pinlabel {$f_{i}|_{D^2_{i}}$} [b] <0pt, 2pt> at 215 521
\pinlabel {$D^2_{i} \subset \R^4$} [bl] at 466 605
\pinlabel {$\pi$} [l] at 419 439
\pinlabel {$l$} [tl] <0pt, 2pt> at 634 5
\pinlabel {$N(l)$} [tl] at 647 -57
\pinlabel* {$\pi \circ f_{i}(D^2_{i})$} [tl] at  357 -120
\pinlabel {$A_{i}$} [b]  at 200 303
\pinlabel {$\R^2$} [bl] at 679 233
\endlabellist
\begin{center}
\begin{minipage}{8cm}
\includegraphics[width=8cm]{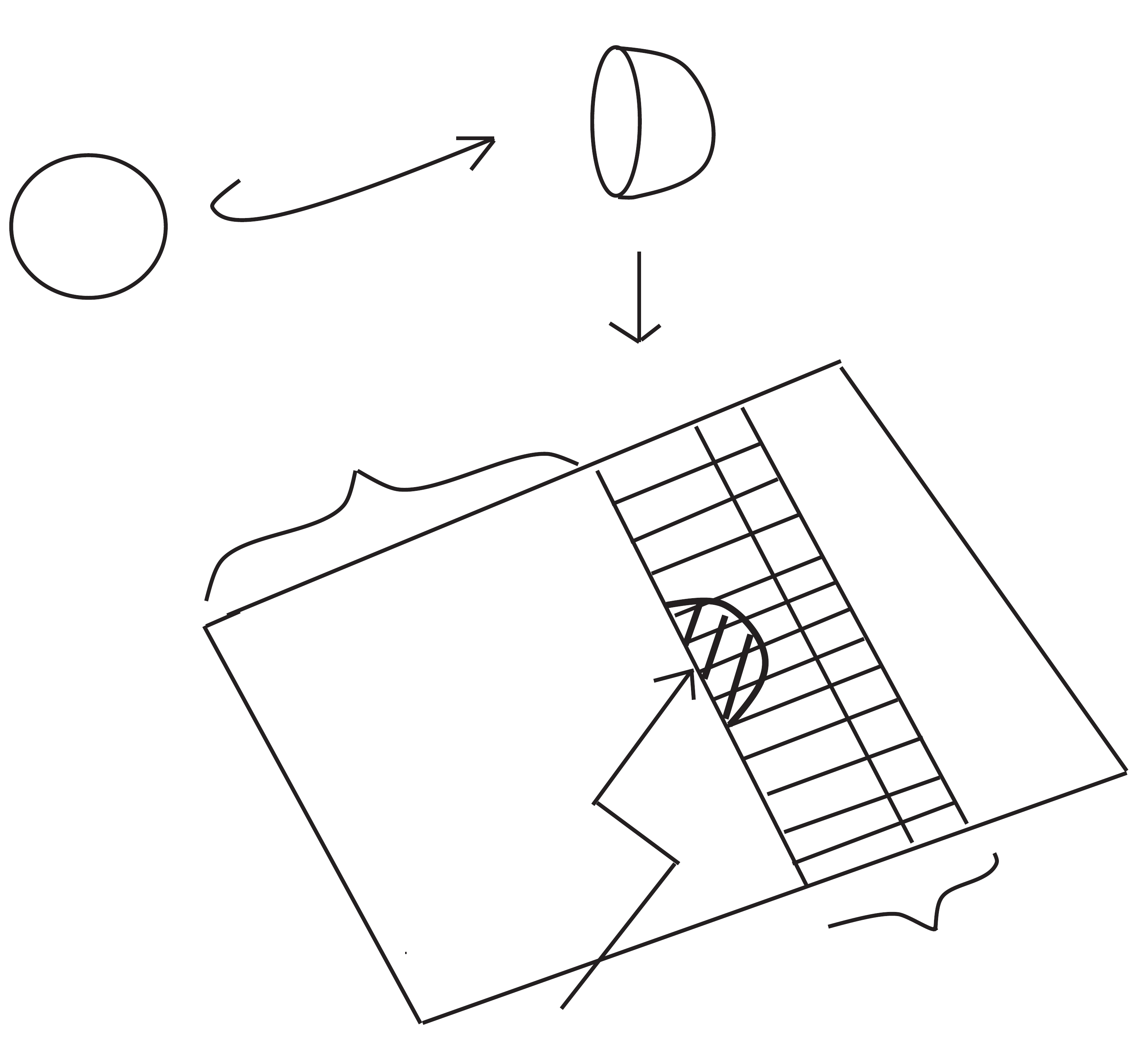}
\end{minipage}
\end{center}
\caption{The mapping $\pi \circ f_{i}|_{D^2_{i}}$}
\label{fig:c}
\end{figure}
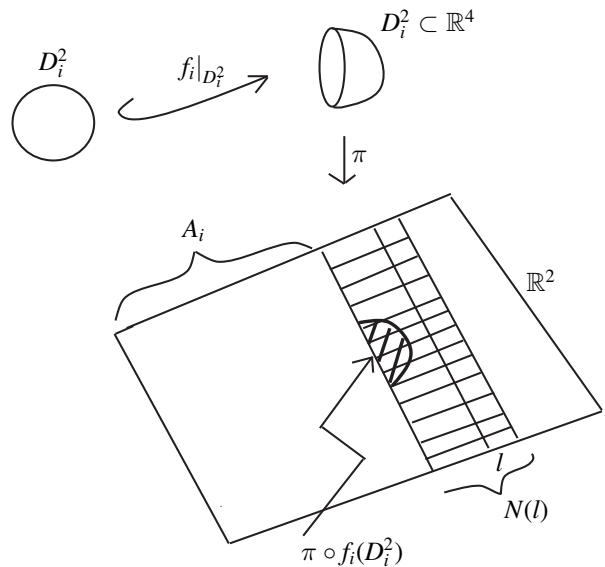

\begin{lemma}\label{lem:t0}
Let $F$ be a closed connected surface and
$f \co F \to \R^4$ be an embedding.
Let $\pi \co \R^4 \to \R^2$ be an orthogonal projection which
is generic with respect to $f$.
Suppose that there exists a proper arc $l$ in $\R^2$
isotopic to a line in $\R^2$ 
such that $\pi \circ f(S(\pi \circ f))$ intersects $l$ transversely at 
two points both of which are the images of fold points.  
Let $N(l)$ be a tubular neighborhood of $l$ in $\R^2$ and
let $A_{0}$ and $A_{1}$ be the connected components of
$\R^2 \setminus $ {\rm Int}$N(l)$.
Then there exist embeddings $f_{i} \co F_{i} \to \R^4$ 
of closed connected surfaces $F_{i}$ into $\R^4$, $i=0,1$,
such that 

\begin{itemize}
\item[{\rm (i)}] $f(F)$ is isotopic to the connected sum 
$f_{0}(F_{0})  \sharp  f_{1}(F_{1})$,
\item[{\rm (ii)}] $\pi$ is generic with respect to $f_{i}$, $i=0,1$,
\item[{\rm (iii)}] $\pi \circ f_{0}(F_{0}) \cap \pi \circ f_{1}(F_{1})
= \emptyset$,
\item[{\rm (iv)}] for $i=0,1$, there exists a $2$--disk $D^2_{i} \subset F_{i}$ 
      such that 
      \begin{itemize}
      \item[{\rm (iv-$1$)}] $F_{i} \setminus 
      {\rm Int}D^2_{i} = (\pi \circ f)^{-1}
      (A_{i})$,
      \item[{\rm (iv-$2$)}] $f_{i}|_{F_{i} \setminus
       {\rm Int}D^2_{i}} = f|_{F_{i}
      \setminus 
      {\rm Int}D^2_{i}}$,
      \end{itemize}
\item[{\rm (v)}] for $i=0,1$,
$\pi \circ f_{i}|_{D^2_{i}}$ is a mapping as depicted in \fullref{fig:c}.                
\end{itemize}

\end{lemma}

\begin{proof}
Set $l_{i} = \partial N(l) \cap A_{i}, i=0,1$. 
Then $(\pi \circ f)^{-1}(l_{i})$ is a closed $1$--dimensional manifold, and the
embedding $f|_{(\pi \circ f)^{-1}(l_{i})}$ into $\pi^{-1}(l_{i})
\cong \R^3$ is a trivial knot, since $\pi \circ f|_{(\pi \circ f)^{-1}(l_{i})}
\co (\pi \circ f)^{-1}(l_{i}) \to l_{i} \cong \R$
is a Morse function with one maximum and one minimum.
Therefore, $f((\pi \circ f)^{-1}(l_{i}))$ bounds a 2--disk $\Delta^2_{i}$
in $\pi^{-1}(l_{i})$, $i=0,1$.
We slightly push the interior of the $2$--disk into $\pi^{-1}({\rm Int}N(l))$
and we denote it by $\tilde{\Delta}^2_{i}$. Then we get the desired embeddings
$f_{i} \co F_{i} = (\pi \circ f)^{-1}(A_{i}) \cup D^2_{i} \to \R^4, i=0,1$,
such that $f_{i}|_{(\pi \circ f)^{-1}(A_{i})} = f|_{(\pi \circ f)^{-1}(A_{i})}$
and $f_{i}(D^2_{i}) = \tilde{\Delta}^2_{i}$. 
\end{proof} 

Let $f \co F \to \R^4$ be an embedding of a closed connected surface $F$ and
$\pi \co \R^4 \to \R^2$ be an orthogonal projection
which is generic with respect to $f$.
Then $\pi \circ f(S(\pi \circ f))$ has
fold crossings and cusps.  
We have four regions locally near a fold crossing, and
we have two regions locally near a cusp.

\begin{lemma}\label{lem:t1}
Let $f \co F \to \R^4$ be an embedding of a closed connected surface $F$
and $\pi \co \R^4 \to \R^2$ be an orthogonal projection which is generic
with respect to $f$.
Then, the local widths around a fold crossing of $\pi \circ f(S(\pi \circ f))$
are of the forms $n, n+2, n+2, n+4$ for some $n \geq 0$ even. 
The local widths around the image of a cusp are of the 
forms $n, n+2$ for some $n \geq 2$ even. See \fullref{fig:d}.

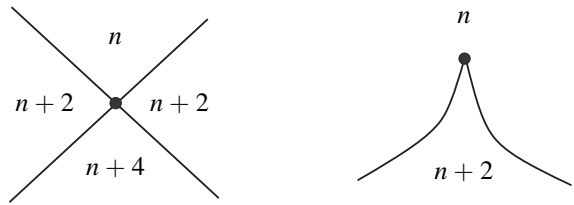
\begin{figure}[ht!]
\begin{center}
\begin{minipage}{4cm}
\labellist\small
\pinlabel $n$ at 328 460
\pinlabel $n+2$ at 213 354
\pinlabel $n+2$ at 432 354
\pinlabel $n+4$ at 328 249
\endlabellist
\includegraphics[width=3cm]{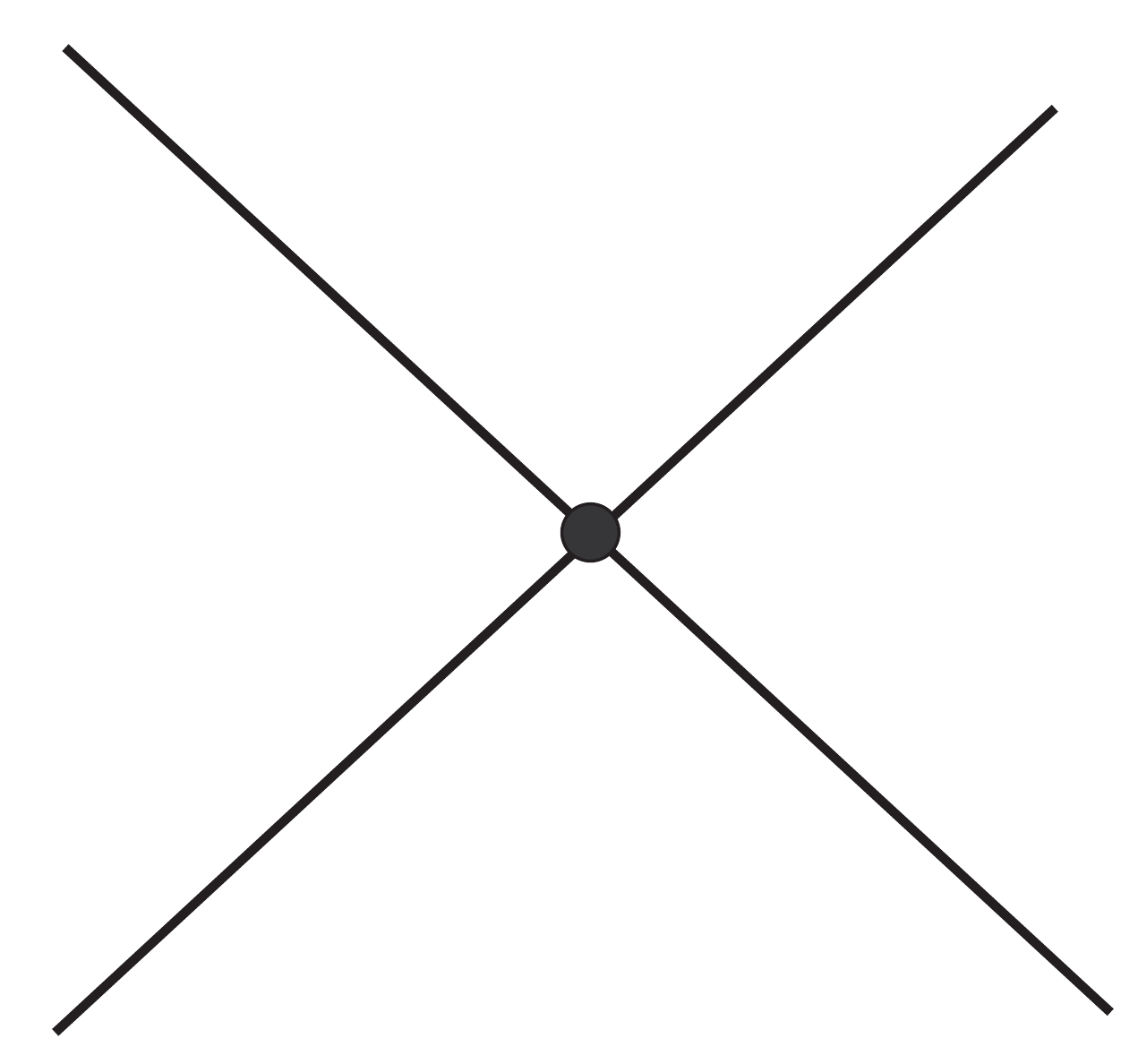}
\end{minipage}
\hspace{0.5cm}
\begin{minipage}{4cm}
\labellist\small
\pinlabel $n$ at 274 502
\pinlabel $n+2$ at 274 249
\endlabellist
\includegraphics[width=3cm]{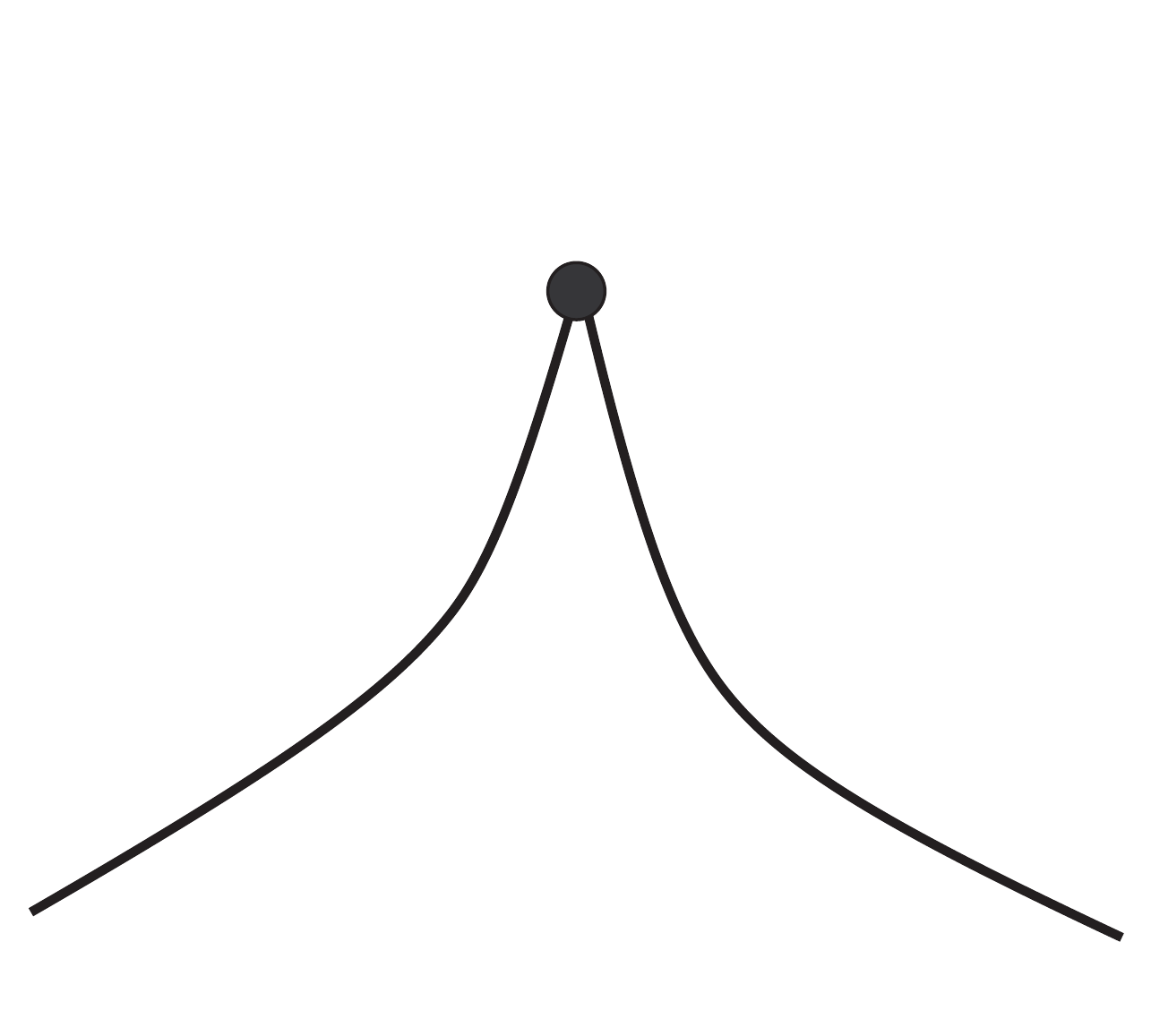}
\end{minipage}
\end{center}
\caption{Local widths around a fold crossing (left) and around a cusp (right)}
\label{fig:d}
\end{figure}
\end{lemma}

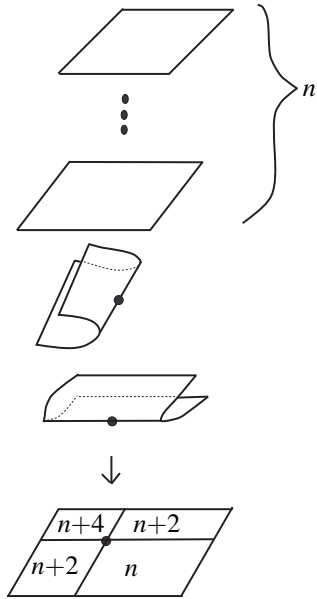
\begin{figure}[ht!]
\begin{center}
\begin{minipage}{4.5cm}
\labellist\small
\pinlabel $n$ [l] at 392 507
\pinlabel $n$ [bl] <0pt,-2pt> at 123 -238
\pinlabel $n{+}2$ [bl] <0pt,-3pt> at 136 -170
\pinlabel $n{+}2$ [bl] <0pt,-2pt> at -20 -238
\pinlabel $n{+}4$ [bl] <0pt,-3pt> at 20 -170
\endlabellist
\includegraphics[width=4.5cm]{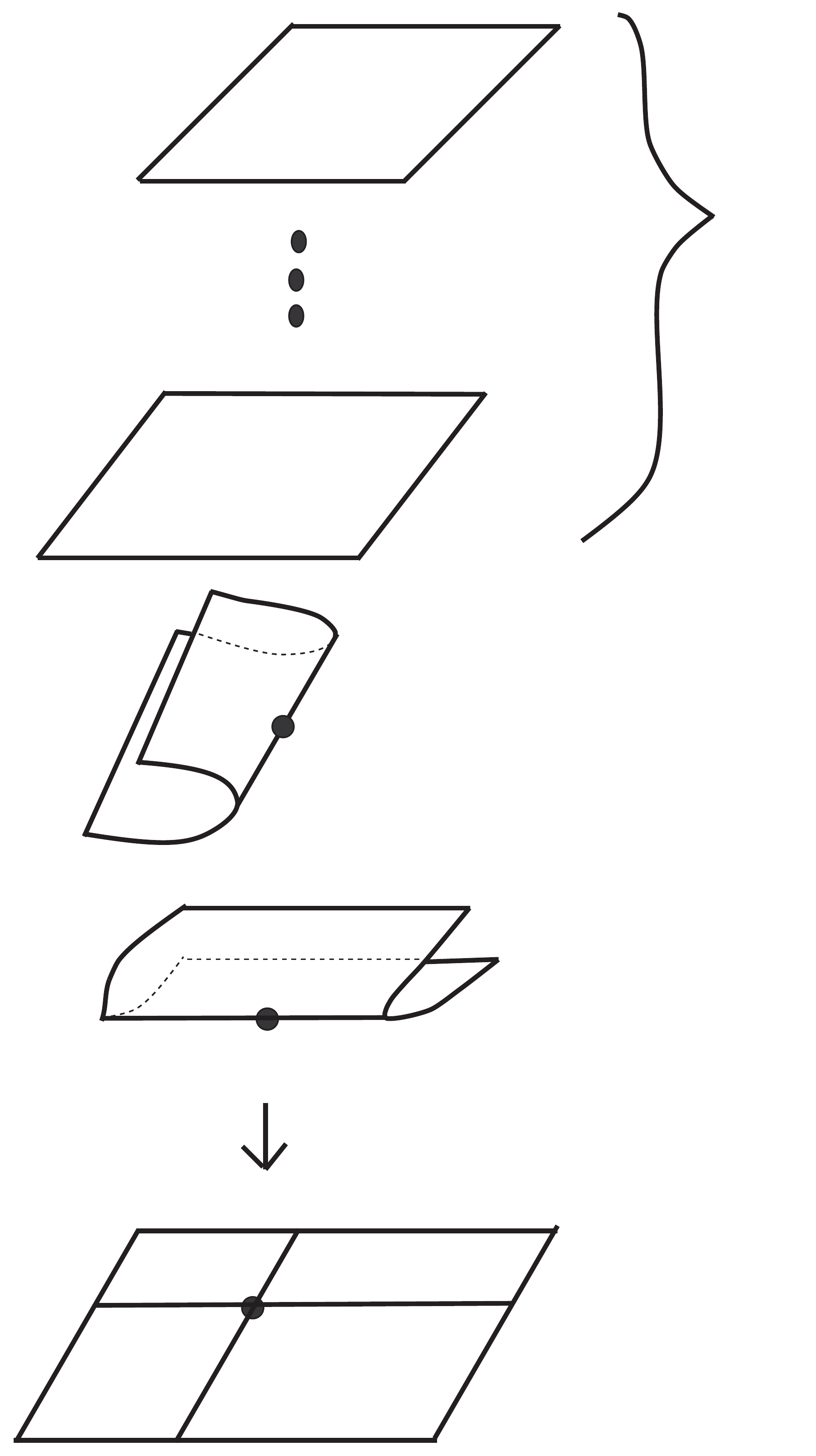}
\end{minipage}
\end{center}
\caption{The situation near a fold crossing}
\label{fig:e}
\vspace{20pt}
\end{figure}

\begin{figure}[ht!]
\begin{center}
\begin{minipage}{4.5cm}
\labellist\small
\pinlabel $n$ [bl] <0pt,-3pt> at 84 -157
\pinlabel $n{+}2$ [bl] <-1pt,-3pt> at 24 -256
\pinlabel $n{-}1$ [l] at 406 464
\endlabellist
\includegraphics[width=4.5cm]{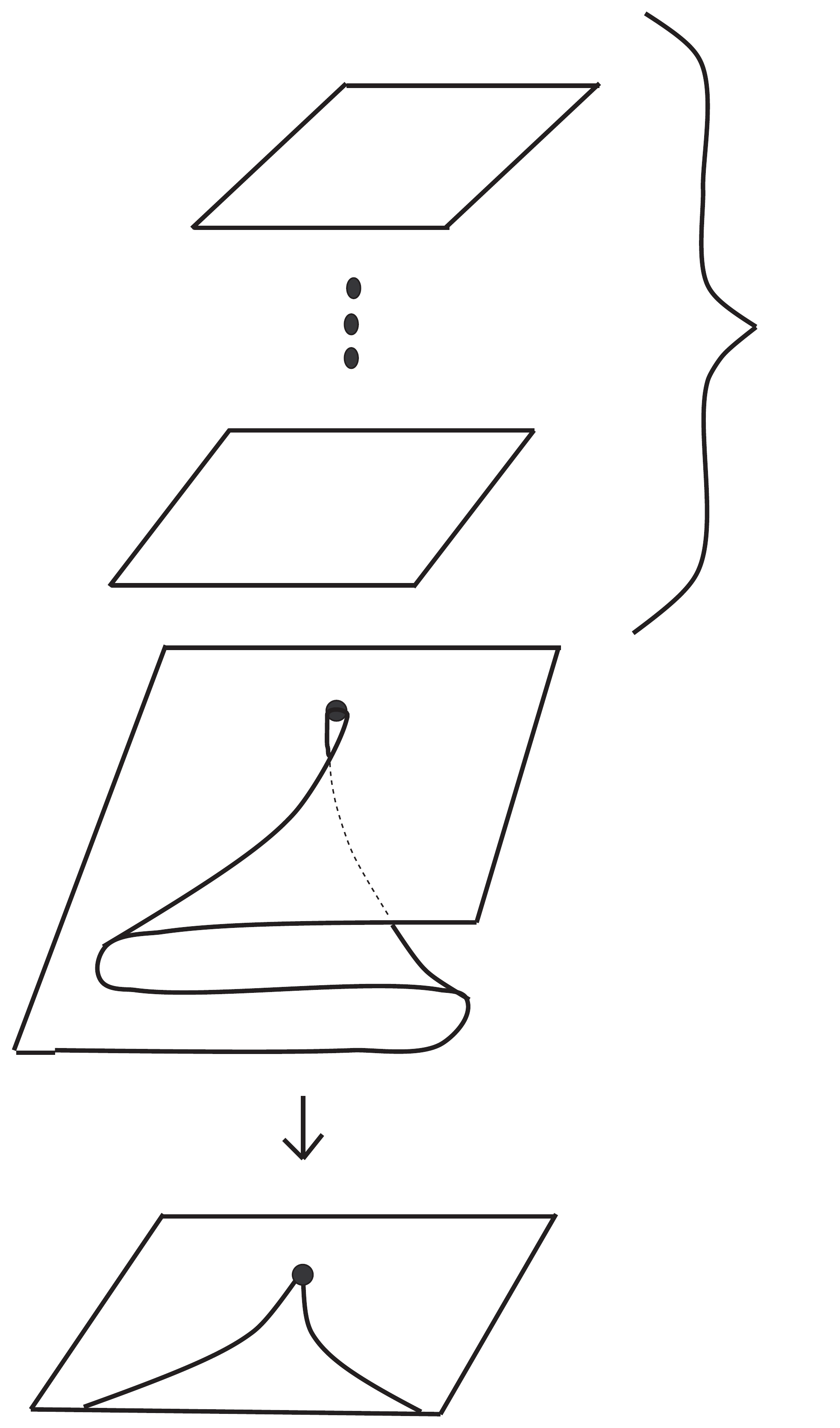}
\end{minipage}
\end{center}
\caption{The situation near a cusp}
\label{fig:f}
\end{figure}

\begin{proof}
If a point $x \in \R^2$ crosses the image of a fold curve,
then the number of elements in the inverse image $(\pi \circ f)^{-1}(x)$
changes by $\pm 2$.
Furthermore, since $F$ is compact, $\pi \circ f$ is not surjective,
and the local width for the unbounded region must be zero.
Therefore, the local width of each region should be an even number.

Let $x \in \R^2$ be a fold crossing.
Then the mapping $\pi \circ f$ near $(\pi \circ f)^{-1}(x)$
is easily seen to be equivalent to the mapping as depicted in \fullref{fig:e}
for some $n \geq 0$.
Furthermore, each local width should be even.
Therefore, the desired conclusion follows.

For a cusp, the situation is as depicted in \fullref{fig:f} for some $n$.
Since the mapping $\pi \circ f$ near a cusp point is an
open map, each local width around the image of a cusp
should be positive.
Then, the desired conclusion follows.
This completes the proof. 
\end{proof}

\vspace{-20pt}
Let us give a characterization of surface knots with width two.

\begin{theorem}\label{thm:t0}
Let $F \subset \R^4$ be a surface knot.
Then $w(F)=2$ if and only if $F$ is strongly trivial.
\end{theorem}

\begin{figure}[ht!]
\psfrag{a}{$2$}
\psfrag{b}{$0$}
\begin{center}
\begin{minipage}{8cm}
\labellist\small
\pinlabel 2 at 238 517
\pinlabel 0 at 238 664
\pinlabel 0 at -155 369
\pinlabel 0 at 132 369
\pinlabel 0 at 378 369 
\pinlabel 0 at 750 369
\endlabellist
\includegraphics[width=8cm]{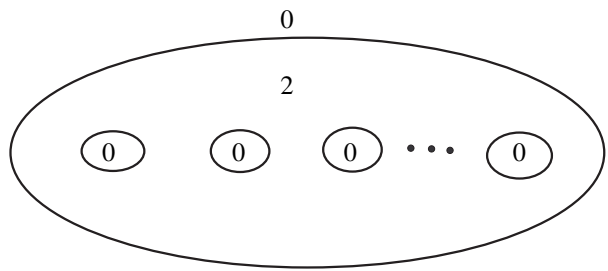}
\end{minipage}
\end{center}
\caption{The image of the singular set of a surface knot $F$
with $w(F)=1$}
\label{fig:g}
\end{figure}

\begin{figure}[ht!]
\begin{center}
\begin{minipage}{8cm}
\labellist\small
\pinlabel 2 at 112 373
\pinlabel 0 at 3 373
\pinlabel 0 at 275 373
\pinlabel $\R$ [b] at 619 192
\endlabellist
\includegraphics[width=8cm]{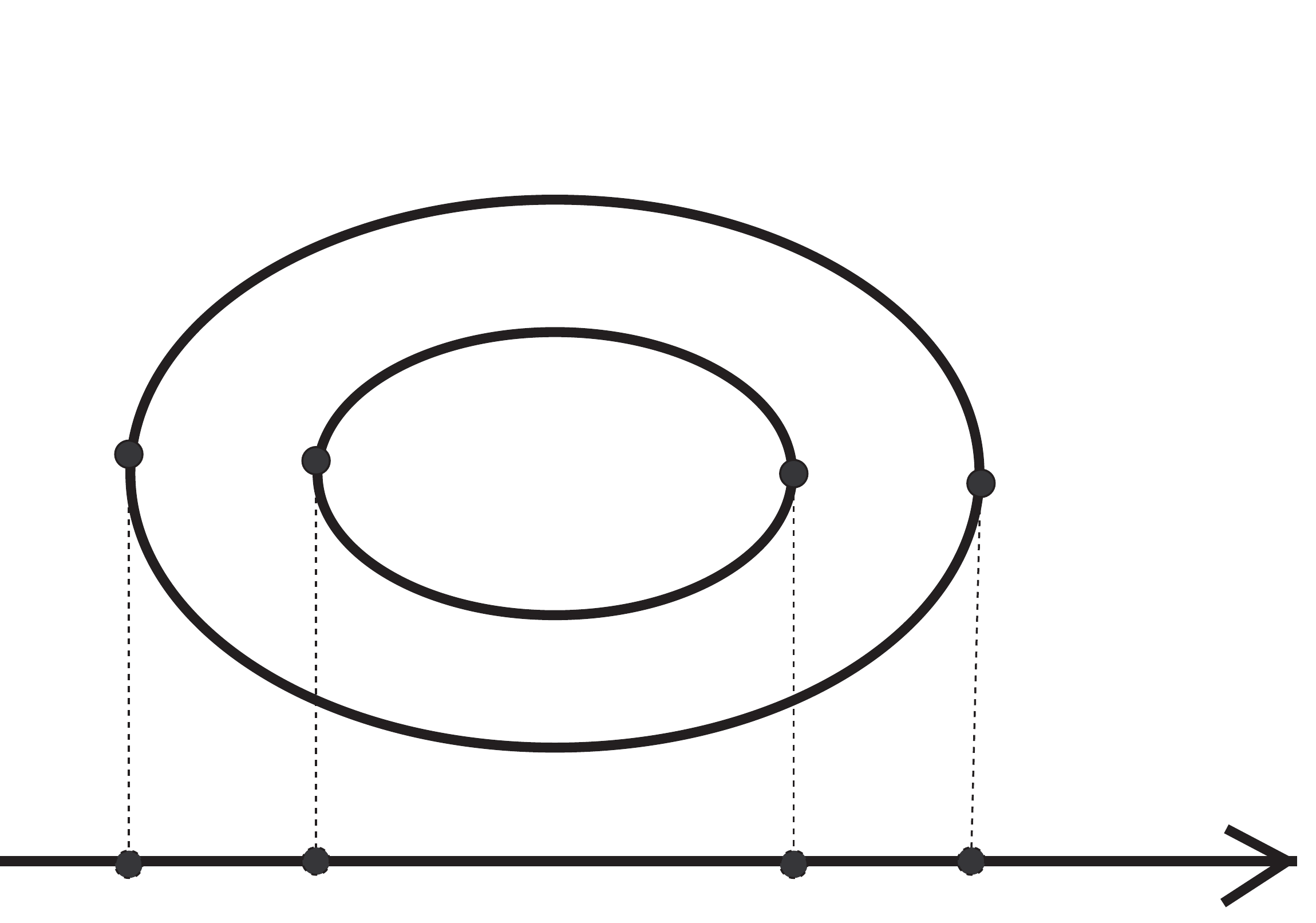}
\end{minipage}
\end{center}
\caption{The image of the singular set $S(\pi|_{F_{j}})$}
\label{fig:h}
\end{figure}

\begin{figure}[ht!]
\begin{center}
\labellist\small
\pinlabel 2 at 232 451
\pinlabel 0 at 21 451 
\pinlabel $\R$ [b] at 460 247
\endlabellist
\includegraphics[width=6cm]{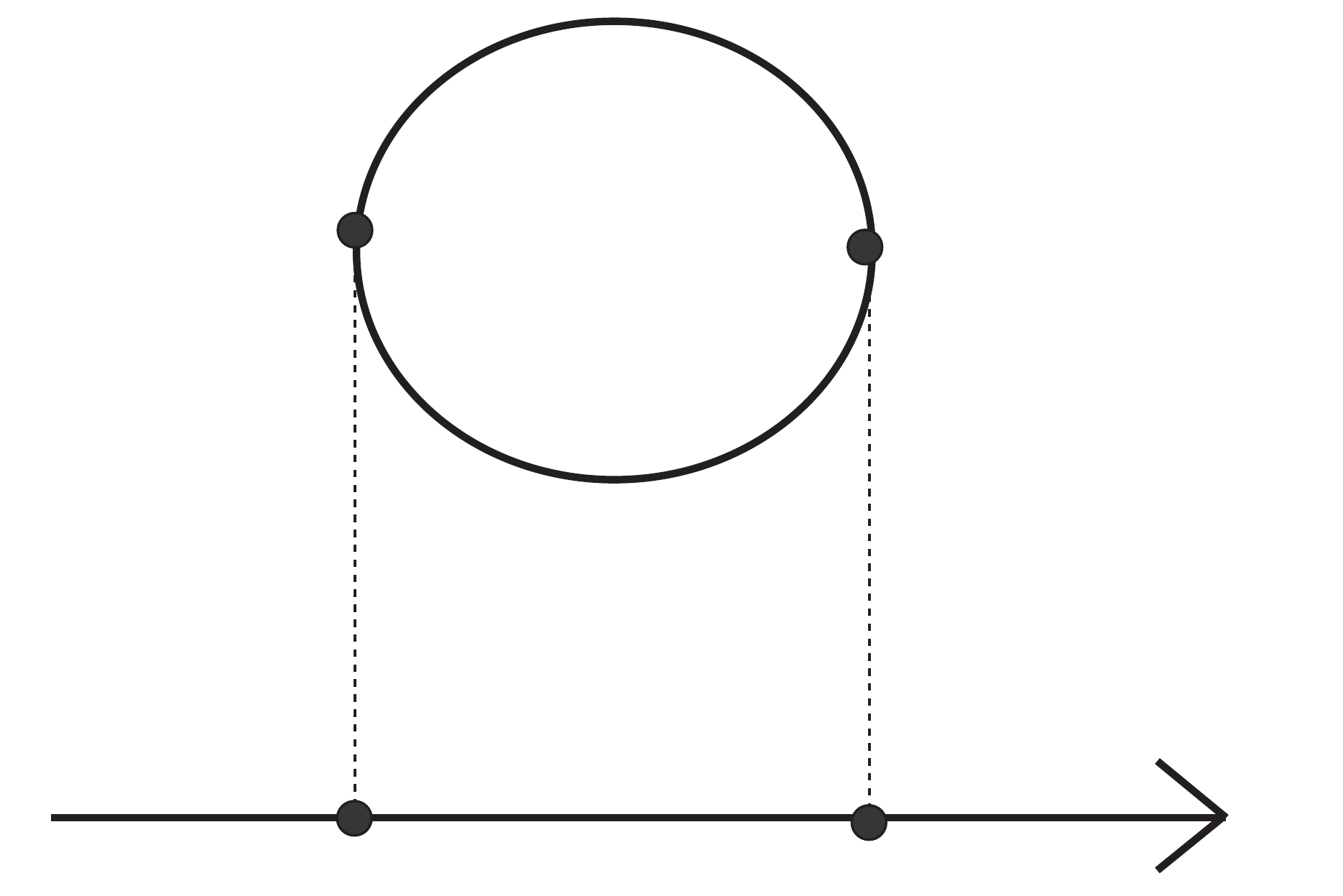}
\end{center}
\caption{The image of the singular set $S(\pi|_{F})$}
\label{fig:i}
\end{figure}

\begin{proof}
Suppose that $w(F)=2$.
We may assume that for an orthogonal projection $\pi \co \R^4 \to \R^2$
which is generic with respect to $F$, we have $w(F,\pi)=w(i,\pi)=2$,
where $i \co F \to \R^4$ is the inclusion mapping.
Then the local width of each region of $\R^2 \setminus \pi(S(\pi \circ i))$
must be equal to $0$ or $2$.
Therefore, by \fullref{lem:t1} there are no fold crossings nor cusps.
Since $F$ is connected, we see that the image of the singular set
$S(\pi|_{F})$ must be as depicted in \fullref{fig:g} up to isotopy of $\R^2$.
Then by using \fullref{lem:t0}, we see that either (i) $F$ is isotopic
to a connected sum $F_{1} \sharp F_{2} \sharp \cdots \sharp F_{r}$
for some $r \geq 1$ such that $\pi$ is generic with respect to $F_{j}$ 
and the image of the singular set $S(\pi|_{F_{j}})$ is as depicted
in \fullref{fig:h} up to isotopy of $\R^2$, $j=1,2, \ldots ,r$,
or (ii) the image of the singular set $S(\pi|_{F})$ is as depicted
in \fullref{fig:i} up to isotopy of $\R^2$.
In case (i), each $F_{j}$ is strongly trivial by \fullref{lem:crs}.
Therefore, $F$ is also strongly trivial.
In case (ii), $F$ is strongly trivial by \fullref{lem:crs}.
Conversely, if $F$ is strongly trivial, then we see easily that
$w(F)=2$.
This completes the proof.
\end{proof}

Let us recall the notion of a ribbon surface knot, which plays an
important role in the theory of surface knots (Cochran, Kamada,
Kawauchi, Tanaka and Yasuda \cite{Co,Ka0,Ka1,Ta1,Ya1}).

\begin{definition}
Let $A=A_{1} \cup A_{2} \cup \cdots \cup A_{k}$
(or $B=B_{1} \cup B_{2} \cup \cdots \cup B_{l}$)
denote a finite disjoint collection of $3$--balls embedded in $\R^4$.
Parametrize each component $B_{i}$ of $B$ as $b_{i} \co D^2 \times [0,1]
\to \R^4$.
Suppose that for each $i=1,2,\ldots,l$, we have
 \begin{itemize}
\item[(i)] $\partial A \cap b_{i}(D^2 \times [0,1])=b_{i}(D^2 \times
\{0,1\})$, and 
\item[(ii)] $b_{i}(D^2 \times (0,1)) \cap A = b_{i}(D^2 \times I_{i})$
for a finite set $I_{i} \subset (0,1)$.
\end{itemize}

Then the surface knot
 \begin{eqnarray}
  F=\left(\partial A \setminus \cup^l_{i=1} b_{i}(D^2 \times 
  \{0,1\})\right)
  \bigcup \cup^l_{i=1}b_{i}(\partial D^2 \times [0,1])    \nonumber 
 \end{eqnarray}
(after a suitable smoothing) is called a \textit{ribbon surface knot}
if $F$ is connected. 
\end{definition}

Note that a surface knot which is strongly trivial is a ribbon surface knot.
If a ribbon surface knot is non-orientable, then
the genus must be even.

\begin{proposition}\label{prop:t0}
Let $F \subset \R^4$ be a ribbon surface knot which is not strongly trivial.
Then we have $w(F)=4$.
\end{proposition}

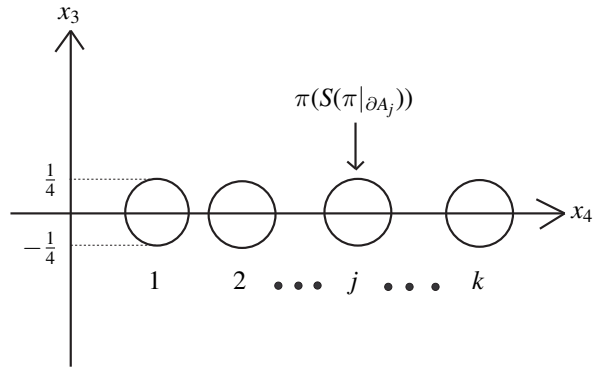
\begin{figure}[ht!]
\begin{center}
\begin{minipage}{8cm}
\labellist\small
\pinlabel $\pi(S(\pi|_{\partial A_{j}}))$ [b] at 336 476
\pinlabel 2 [B] at 179 250
\pinlabel 1 [B] at 62 250
\pinlabel $\frac{1}{4}$ [r] at -60 401
\pinlabel ${-}\frac{1}{4}$ [r] at -60 307
\pinlabel $j$ [B] at 336 250
\pinlabel $k$ [B] at 507 250
\pinlabel $x_{3}$ at -55 629
\pinlabel $x_{4}$ [l] at 624 354
\endlabellist
\includegraphics[width=8cm]{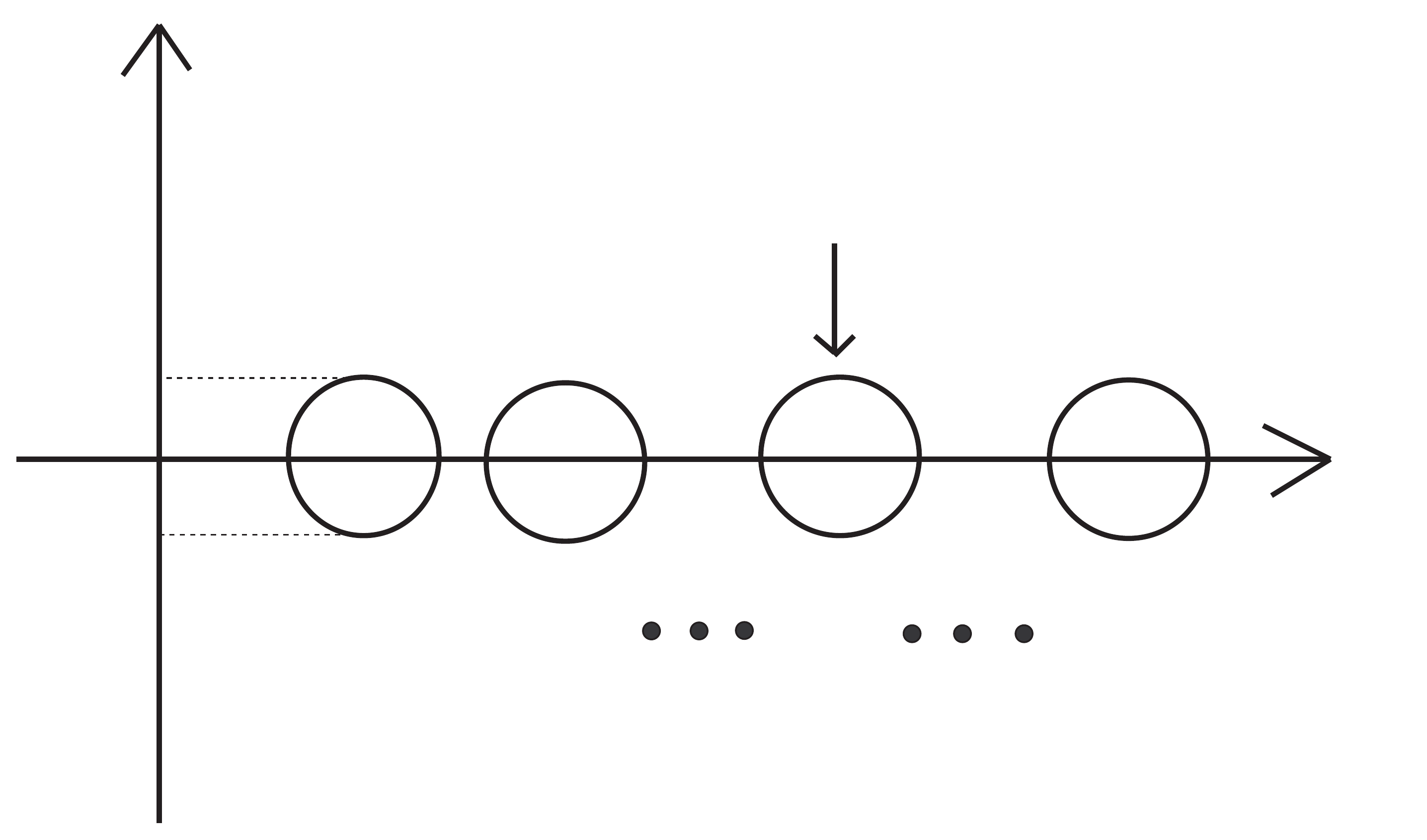}
\end{minipage}
\end{center}
\caption{$\pi(S(\pi|_{\partial A_{j}}))$} 
\label{fig:j}
\end{figure}

\begin{figure}[ht!]
\begin{center}
\begin{minipage}{5cm}
\labellist\small
\pinlabel $\pi(A)$ [bl] at 432 560
\pinlabel {$\pi \circ b$ restricted to ${\displaystyle \amalg^l_{i=1}
(\{0\} \times [0,1])_{i}}$} [b] at 610 883
\endlabellist
\includegraphics[width=5cm]{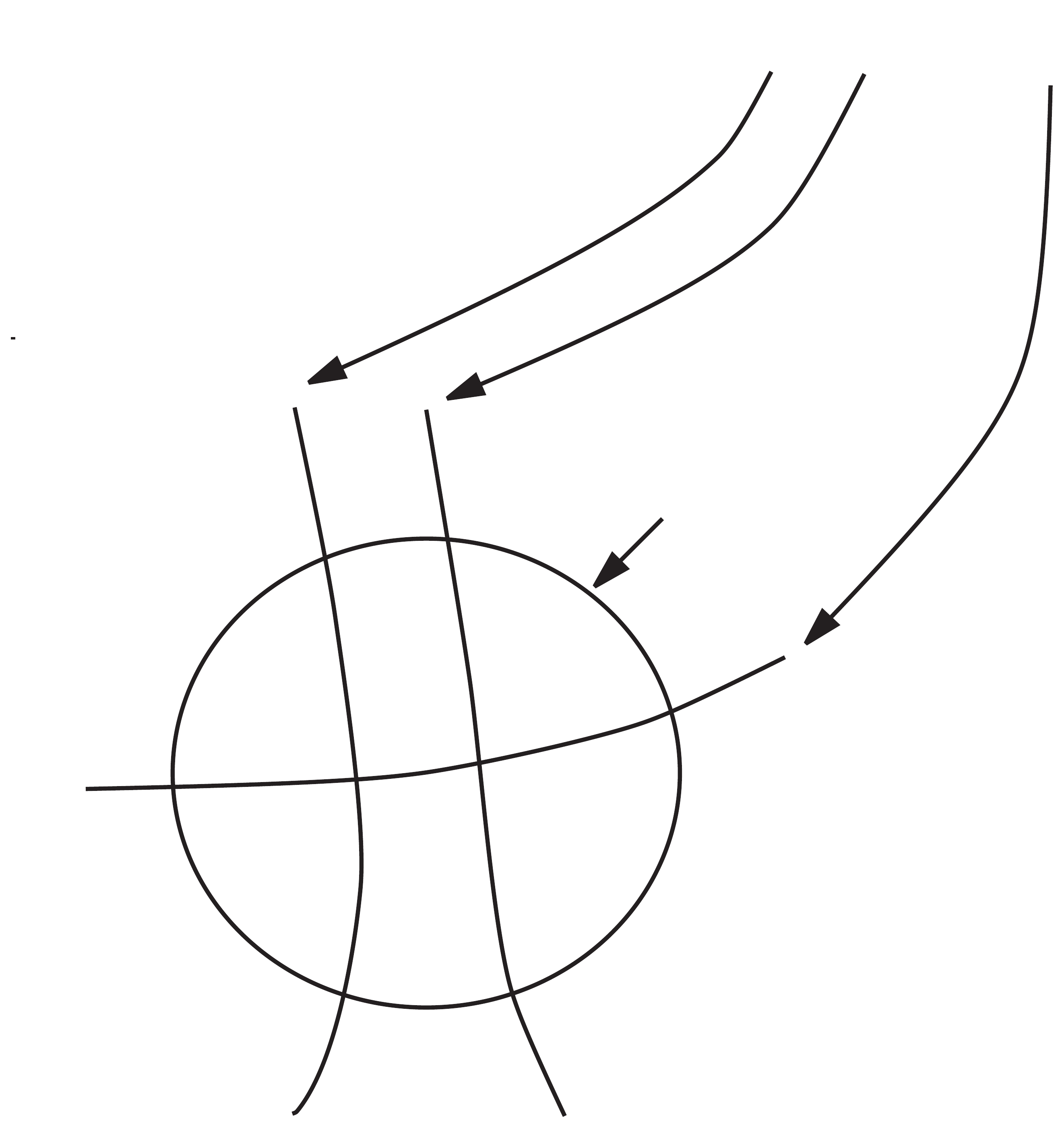}
\end{minipage}
{\Huge $\Rightarrow$}
\begin{minipage}{5cm}
\includegraphics[width=5cm]{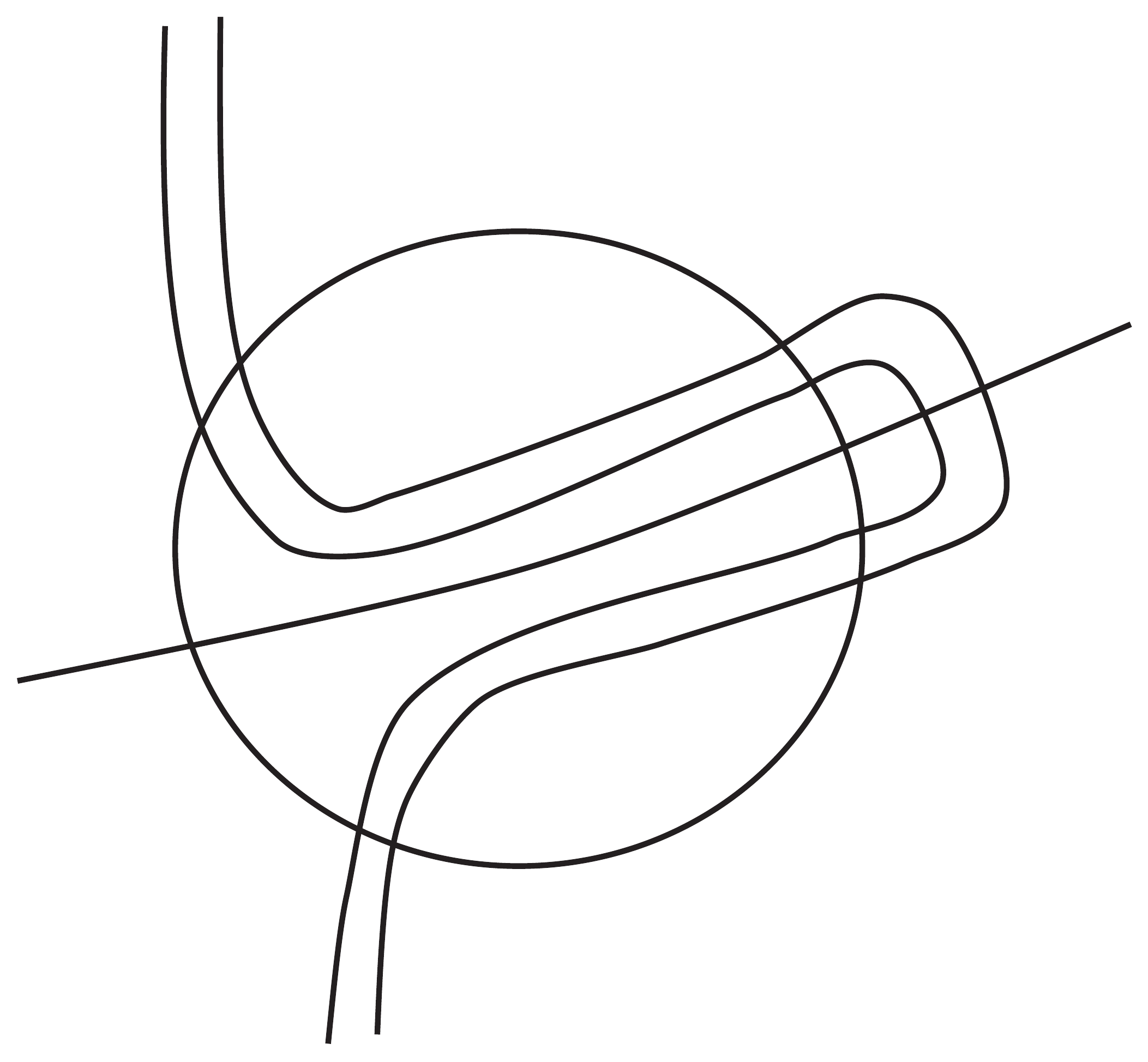}
\end{minipage}
\end{center}
\caption{Pushing the crossings of $\pi \circ b$ restricted to ${\displaystyle 
\amalg^l_{i=1}(\{0\} \times [0,1])_{i}}$ out of $\pi(A)$} 
\label{fig:k}
\vspace{20pt}
\end{figure}

\begin{figure}[ht!]
\begin{center}
\begin{minipage}{8cm}
\includegraphics[width=8cm]{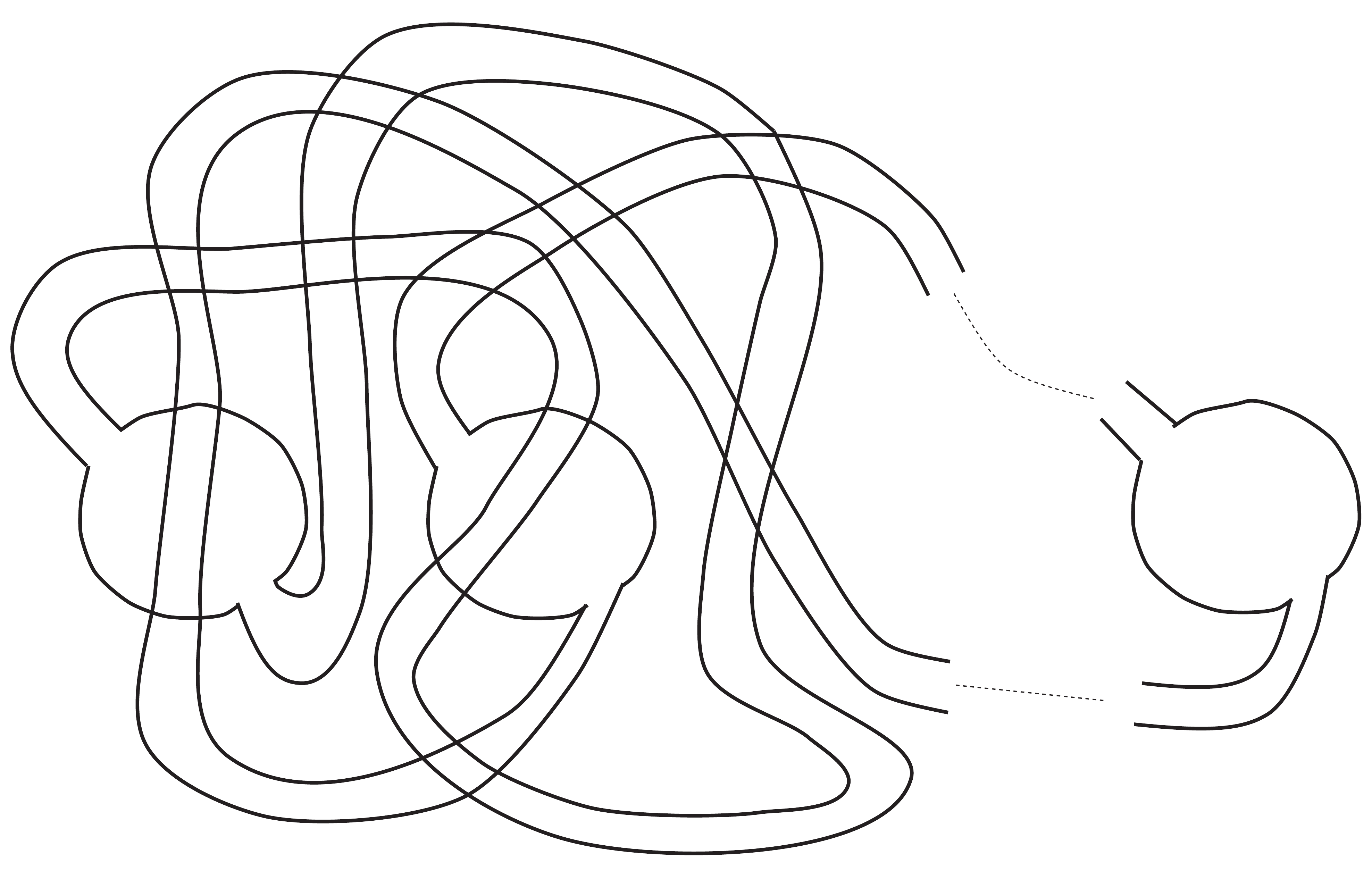}
\end{minipage}
\end{center}
\caption{$\pi(S(\pi|_{F}))$} 
\label{fig:l}
\vspace{20pt}
\end{figure}

\begin{proof}
By isotopy of $F$ we may assume that 
$$A_{j}=\{(x_{1},x_{2},x_{3},x_{4}) \in \R^4 \ | \ x_{1}=0, x^2_{2}+x^2_{3}+(x_{4}-j)^2 \leq (1/4)^2 \}, \ j=1,2,\ldots,k.$$
We define $\pi \co \R^4 \to \R^2$ by $\pi(x_{1},x_{2},x_{3},x_{4})
=(x_{3},x_{4})$.
Then $\pi$ is generic for $\partial A_{j}$ and $\pi(S(\pi|_{\partial  A_{j}}))$
is as depicted in \fullref{fig:j}. 
Moreover, we may further assume
that each $b_{i}|_{\{0\} \times [0,1]}$ satisfies
$b_{i}(0,0), b_{i}(0,1) \in S(\pi|_{\partial A})$ and 
$b_{i}|_{\{0\} \times I_{\varepsilon }}$ 
is an embedding into the closure of 
$\{(0,0,x_{3},x_{4})\in \R^4\} \setminus A$, where $I_{\varepsilon}
=[0,\varepsilon) \cup (1-\varepsilon,1]$ and
$\varepsilon > 0$ is sufficiently small.

We define 
$$b \co \displaystyle \amalg^l_{i=1}(D^2 \times [0,1])_{i} \to \R^4$$
by $b(x)=b_{i}(x), x \in (D^2 \times [0,1])_{i}$, where 
$(D^2 \times [0,1])_{i}$ is a copy of $D^2 \times [0,1]$, $i=1,2,\ldots,l$.
\vspace{5pt}

We may assume that
$\pi \circ b$ restricted to ${\displaystyle \amalg^l_{i=1}(\{0\} \times [0,1]
)_{i}}$ is an 
immersion with normal crossings.
Furthermore, by pushing the crossings out of $\pi(A)$ one by one
by an isotopy of $F$, we may assume that $\pi(A)$ does not contain 
any double point of 
$\pi \circ b$ restricted to 
${\displaystyle \amalg^l_{i=1}(\{0\} \times [0,1])_{i}}$
(see \fullref{fig:k}).

Now the fiber of the normal disk bundle to 
$b_{i}(\{0\} \times [0,1])$ in $\R^4$
is a $3$--dimensional disk.
If we fix $b_{i}(\{0\} \times [0,1])$, then the
isotopy class of $b_{i}(D^2 \times [0,1])$ is determined by
the homotopy class of a unit normal vector field along
$b_{i}(\{0\} \times [0,1])$, which corresponds to the unit normal vector to
$b_{i}(D^2 \times \{*\})$ in the $3$--dimensional disk fiber.
Therefore, we may assume that the tangent plane to 
$b_{i}(D^2 \times \{t\})$ at $b_{i}(\{0\} \times \{t\})$
is not parallel to the fibers of $\pi \co \R^4 \to \R^2$,
$t \in [0,1]$.

By taking 
$B={\displaystyle \amalg^l_{i=1}b_{i}(D^2 \times [0,1])}$
``thin" enough, we may then assume that 
$S(\pi \circ b_{i}|_{\partial D^2 \times [0,1]})$
consists exactly of two arcs for each $i$.
Now $\pi(S(\pi|F))$ is as depicted in \fullref{fig:l} and 
we see that the local width of each region is equal to
$0,2$ or $4$.
Therefore, we have $w(F) \leq 4$.
Then by \fullref{thm:t0}, the desired conclusion follows.
This completes the proof.
\end{proof}

Let us recall the notion of bridge index for classical knots.
Here, we give a definition suitable for our purpose.

\begin{definition}
Let $K$ be a classical knot and $\pi \co \R^3 \to \R$ a generic
orthogonal projection.
Let $m(K,\pi)$ be the number of local maxima of $\pi|_{K} \co K \to \R$.
Then the \textit{bridge index} $b(K)$ of $K$ is 
defined to be the minimum of $m(\tilde{K},\tilde{\pi})$,
where $\tilde{K}$ runs through all embeddings of $S^1$ into $\R^3$
isotopic to $K$, and $\tilde{\pi}$ runs through all orthogonal projections
$\R^3 \to \R$ generic with respect to $\tilde{K}$.
A knot having bridge index $n$ is called an \textit{$n$--bridge knot}.
\end{definition}

Note that an orthogonal projection $\pi \co \R^3 \to \R$ is {\it generic
with respect to $K$} if $\pi|_{K} \co K \to \R$ has only non-degenerate 
critical points as its singularities.

\begin{figure}[ht!]
\begin{center}
\labellist\small
\pinlabel {$\R^3_{+}$} <2pt,0pt> at 530 348
\pinlabel {$\R^2$} [b] at 28 785
\pinlabel {once} [l] <0pt,3pt> at 73 733
\pinlabel {$n$ times} [l] <0pt,3pt> at 363 566
\pinlabel* {$k$} [tl] at 334 19
\endlabellist
\includegraphics[width=5cm]{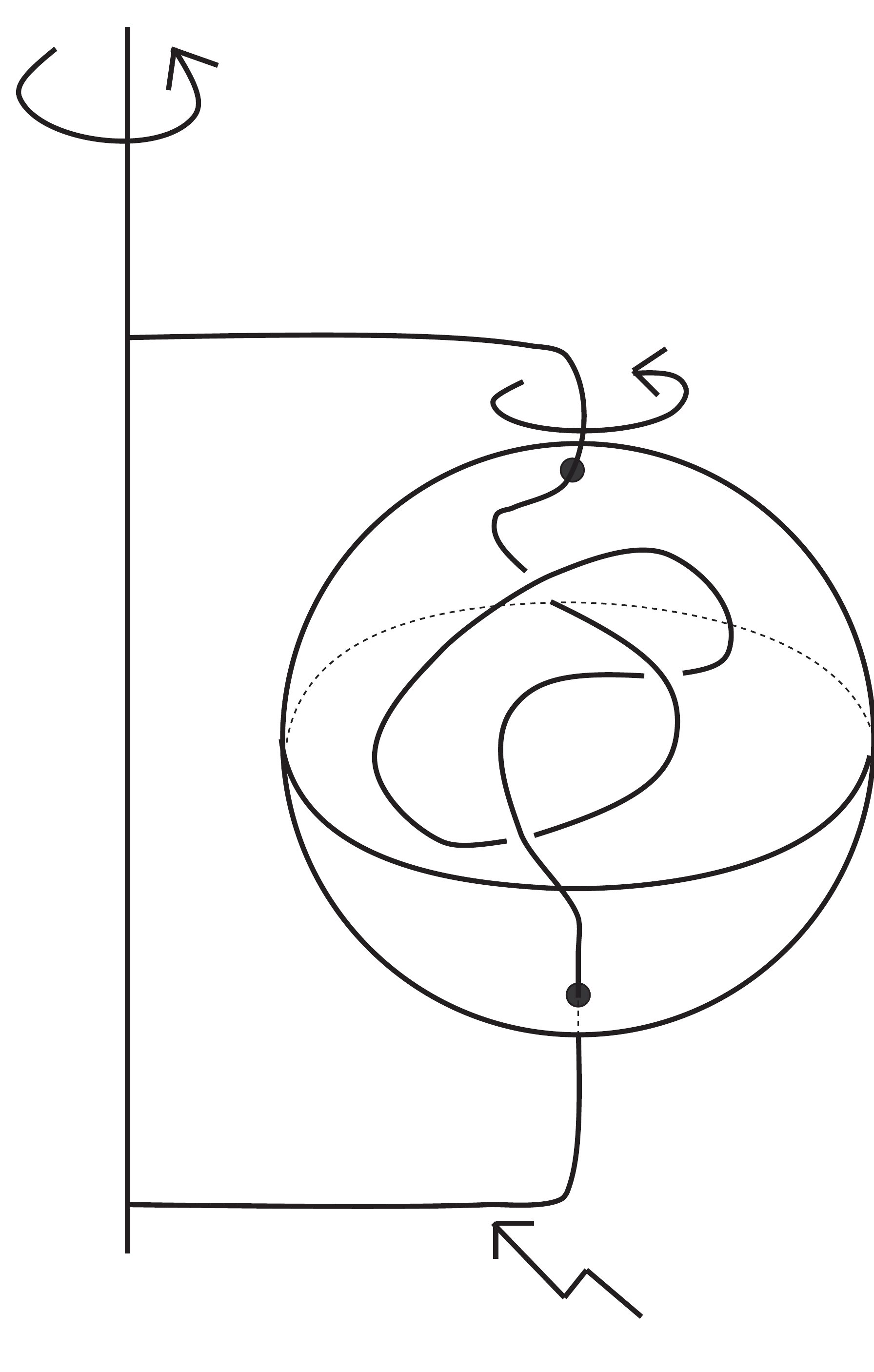}
\end{center}
\caption{The $n$--twist spun trefoil} 
\label{fig:m}
\vspace{20pt}
\end{figure}
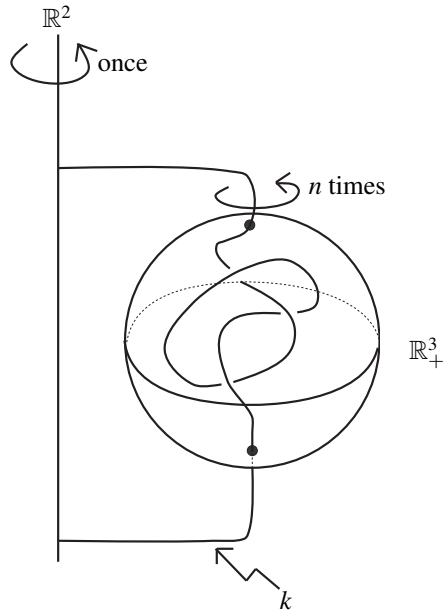

\begin{definition}
Let $\R^3_{+}$ be the $3$--dimensional upper half-space,
ie, $$\R^3_{+}=\{(x_{1},x_{2},x_{3},x_{4})
 \ | \ x_{3} \geq 0, x_{4}=0\}$$ and $\R^2$ the plane $\R^2=\{(x_{1},x_{2},x_{3},x_{4}) \ | \ x_{3}=0, x_{4}=0\}$.
Let $k$ be an arc properly embedded in the half-space 
$\R^3_{+}$.
When the half-space is rotated around the plane 
$\R^2$ in $\R^4$,
the continuous trace of $k$ forms a $2$--sphere.
This $2$--sphere is said to be derived from $k$ by (untwisted) spinning,
and we call the resulting surface knot a \textit{spun knot}.
Moreover, put the knotted part of $k$ in a $3$--ball as in \fullref{fig:m} and
twist it $n$ times, $n \in \Z$, as the half-space spins once around $\R^2$.
Then we call the resulting surface knot an \textit{$n$--twist spun knot}.
In general, $k$ is associated with a knot in $\R^3$,
which is obtained by connecting
the end points of $k$ in an obvious way by an arc in $\R^2$.
See also Zeeman \cite{Ze}.
\end{definition}

\begin{proposition}\label{prop:t1}
Let $F \subset \R^4$ be an $n$--twist spun $2$--bridge knot
with $n \neq \pm 1$.
Then we have $w(F)=4$.
\end{proposition}

\begin{figure}[ht!]
\begin{center}
\labellist\small
\pinlabel {$\R$} [b] at -185 723 
\pinlabel {$a_{3}$} [r] at -185 650 
\pinlabel {$a_{2}$} [r] at -185 510
\pinlabel {$\pi$} [b] at -105 380
\pinlabel {$a_{1}$} [r] at -185 143 
\pinlabel {$0$} [r]  <-4pt, 1pt> at -185 30 
\pinlabel {$a_{0}$} [r] at -185 -90
\pinlabel {$K^{'}$} [l] at 238 592
\pinlabel {$3$--string braid} [bl]  <-4pt, 0pt> at 441 611
\pinlabel {$x_{1}=x_{2}=0$} [l] <0pt, 1.5pt> at 588 516
\pinlabel {$x_{3}=0$} [l] <0pt, 1.5pt> at 748 63
\pinlabel {``knotted part"} [b] at 189 766
\endlabellist
\includegraphics[width=8cm]{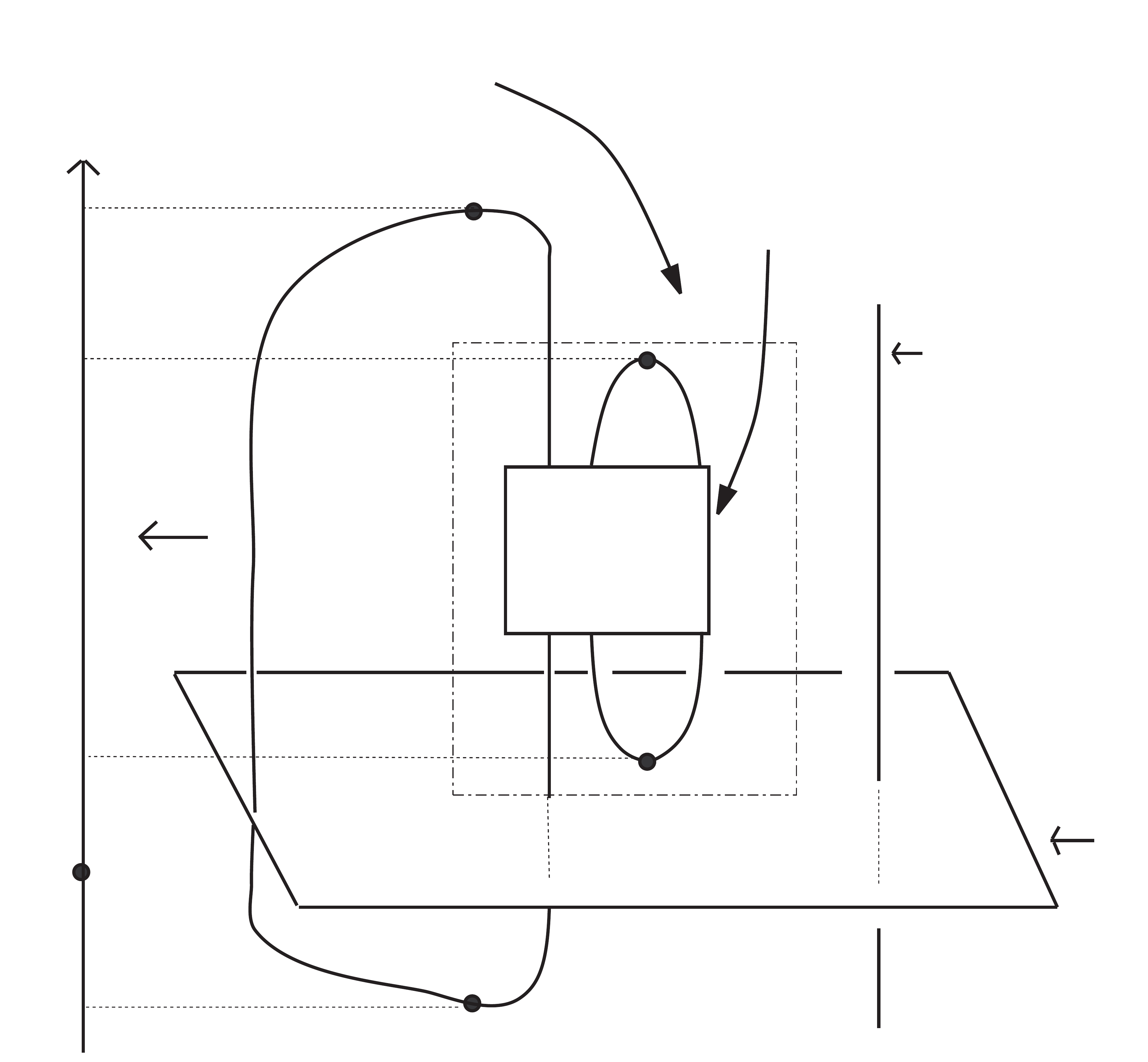}
\end{center}
\caption{Bridge presentation} 
\label{fig:n}
\end{figure}
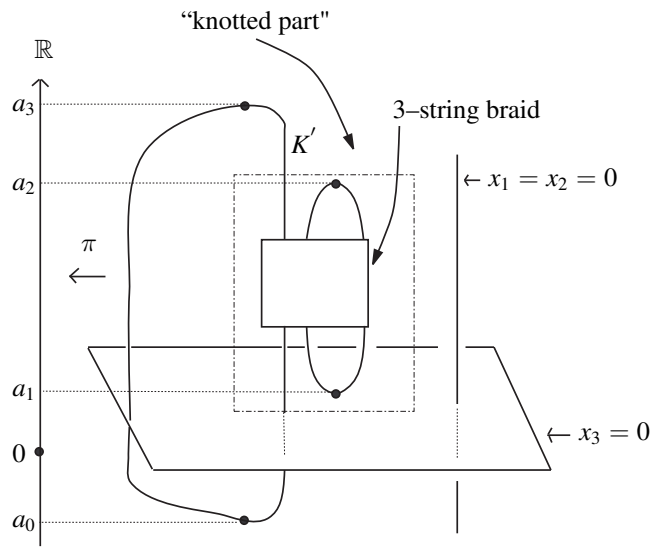

\begin{figure}[ht!]
\begin{center}
\labellist\small
\pinlabel {$0$} at 516 613
\pinlabel {$2$} at 438 543
\pinlabel {$2$} at 326 391
\pinlabel {$4$} at 391 473
\pinlabel {$x_{3}$} [l] at 742 354
\pinlabel {$x_{4}$} [b] at 281 738 
\pinlabel* {$a_{1}$} [tl] at 453 92
\pinlabel* {$a_{2}$} [tl] at 516 129
\pinlabel* {$a_{3}$} [tl] at 584 135
\endlabellist
\includegraphics[width=6cm]{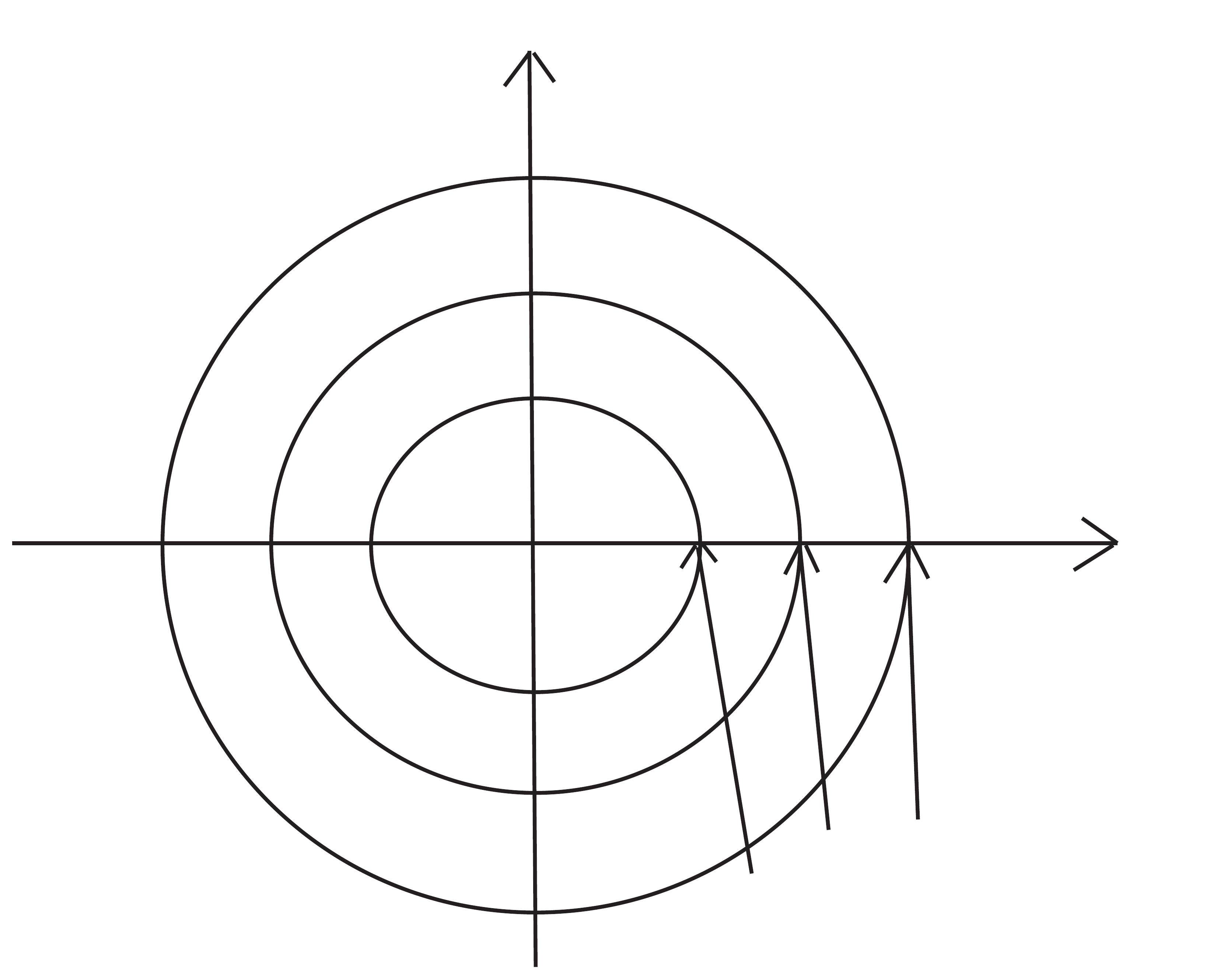}
\end{center}
\caption{A planar projection of a $0$--twist spun $2$--bridge knot and the
associated local widths} 
\label{fig:o}
\end{figure}

\begin{proof}
Let $K$ be a $2$--bridge knot and $\pi \co \R^3 \to \R$ the 
orthogonal projection defined by $\pi(x_{1},x_{2},x_{3}) = x_{3}$.
Then there exists a knot $K^{'}$ isotopic to $K$ such that
$\pi$ is generic for $K^{'}$ and $\pi|_{K^{'}}$ has
two local minima $a_{0}, a_{1}$ and two local maxima $a_{2}, a_{3}$
with $a_{0} < 0 < a_{1} < a_{2} < a_{3}$.
We may assume that the values of the local maxima and the
local minima are all distinct and that
$K^{'}$ is in a position as described in \fullref{fig:n}. 
Rotate the part $K^{'} \cap \R^3_{+} = K^{'} \cap \{x_{3} \geq 0, x_{4}=0 \}$
around $\R^2=\{x_{3}=x_{4}=0\}$ in $\R^4$.
Then we get the $0$--twist spun $F_{0}$ of $K$.
The orthogonal projection $\tilde{\pi} \co \R^4 \to \R^2$ 
defined by $\tilde{\pi}(x_{1},x_{2},x_{3},x_{4})
=(x_{3},x_{4})$ is generic for $F_{0}$ and $\tilde{\pi}(S(
\tilde{\pi}|_{F_{0}}))$ is as depicted in \fullref{fig:o}.
Therefore, we have $w(F_{0}) \leq 4$.
For the $n$--twist spun $F_{n}$ of $K$, rotate $K^{'} \cap \R^3_{+}$
around $\R^2$ once and twist the ``knotted part" $n$ times.
Then $\pi|_{K^{'}}$ does not change and the image of the singular
set is again as depicted in \fullref{fig:o}.
Therefore, we have $w(F_{n}) \leq 4$.
If $w(F_{n})=2$, then by \fullref{thm:t0} $F_{n}$ is strongly trivial.
However, for $n \neq \pm 1$, it is known that $F_{n}$
is not strongly trivial (Cochran \cite{Co}).
Therefore, we have $w(F_{n})=4$ for $n \neq  \pm 1$.
\end{proof}

\begin{remark}
By an argument similar to that in the proof of \fullref{prop:t1},
we can show that the width of
an $n$--twist spun $m$--bridge knot is smaller than or equal to $2m$.
However, even if $n \neq \pm 1$, the equality may not hold.
In fact, for every knot,
its $0$--twist spun is a ribbon surface knot (see, for example, \cite{Co}).
Hence, by \fullref{prop:t0}, we have $4=w(F) < 2m$ if $F$
is a $0$--twist spun $m$--bridge knot with $m \geq 3$.
\end{remark}

\section{Braid index and width}\label{sec-4}

In this section, we study the relationship between 
the braid index and the width of a surface knot.
Throughout this section, we assume that surface knots are orientable.

The notion of surface braid was introduced by Kamada \cite{Ka}.
Kamada and Viro showed that every orientable surface knot
is isotopic to a simple closed surface braid.

A {\it closed surface braid} in $D^2 \times S^2 \subset
(D^2 \times S^2) \cup (D^3 \times S^1)=S^4$ is a closed oriented surface $F$
embedded in $D^2 \times S^2$ such that the restriction map
$pr_{2}|_{F} \co F \to S^2$ of the projection 
$pr_{2} \co D^2 \times S^2 \to S^2$ to the second factor 
is an orientation preserving branched covering.
We say that it is a {\it simple} closed surface braid if $pr_{2}|_{F}$
is a simple branched covering.
An orientation preserving branched covering $f \co F \to M$
between closed oriented
surfaces is {\it simple} if for every branch point $y \in M$,
we have $\sharp f^{-1}(y) = {\rm deg}(f)-1$, where $\sharp$ denotes the
number of elements and deg($f$)$> 0$ is the mapping degree of $f$.
The mapping degree of $pr_{2}|_{F} \co F \to S^2$ is called the 
{\it degree} of the closed surface braid.

The {\it braid index} Braid($F$) of an oriented surface knot $F$ in $\R^4$
is the minimum degree of simple closed surface braids in 
$S^4=\R^4 \cup \{\infty\}$ that
are isotopic to $F$.

For classical knots, the bridge index is smaller than or equal to the 
braid index.
On the other hand, the relation  between the width and the braid
index for classical knots has not been studied as far as the author knows.

For surface knots, we have the following.
 
\begin{proposition}
Let $F \subset \R^4$ be an orientable surface knot.
Then we have
\begin{eqnarray}
 w(F) \leq 2({\rm Braid}(F) + 1).   \nonumber
\end{eqnarray}
\end{proposition}

\begin{figure}[ht!]
\begin{center}
\labellist\small
\pinlabel {image of the singular points of $pr_{2}|_{F}$} [b] at 248 800
\pinlabel {image of the fold curves} [t] at -300 170
\endlabellist
\includegraphics[width=6cm]{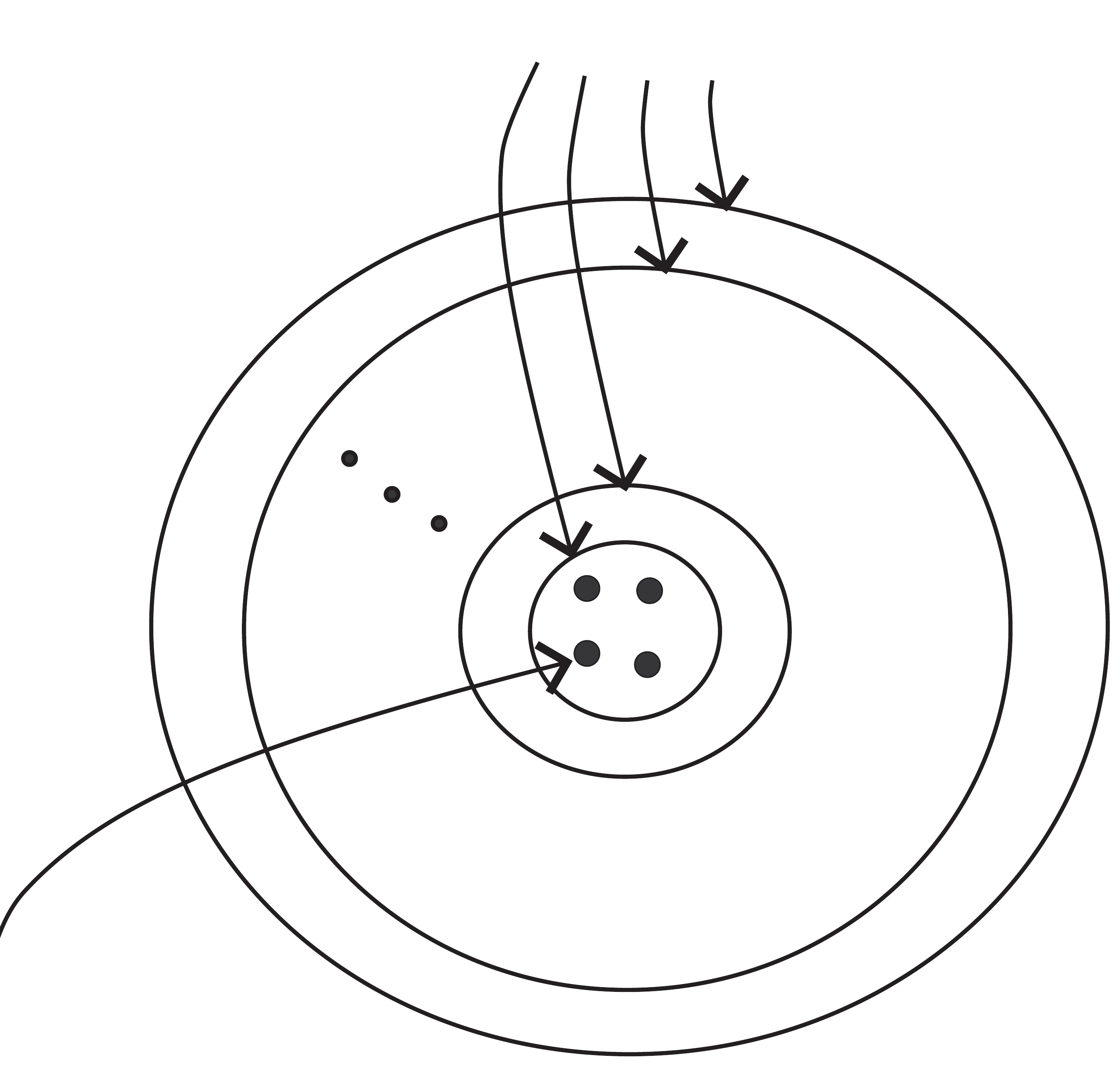}
\end{center}
\caption{Image of the singular points of $\pi|_{F}$} 
\label{fig:p}
\vspace{10pt}
\end{figure}

\begin{figure}[ht!]
\begin{center}
\noindent\llap{\raise 40pt\hbox{(1)}\qquad}\includegraphics[width=6cm]{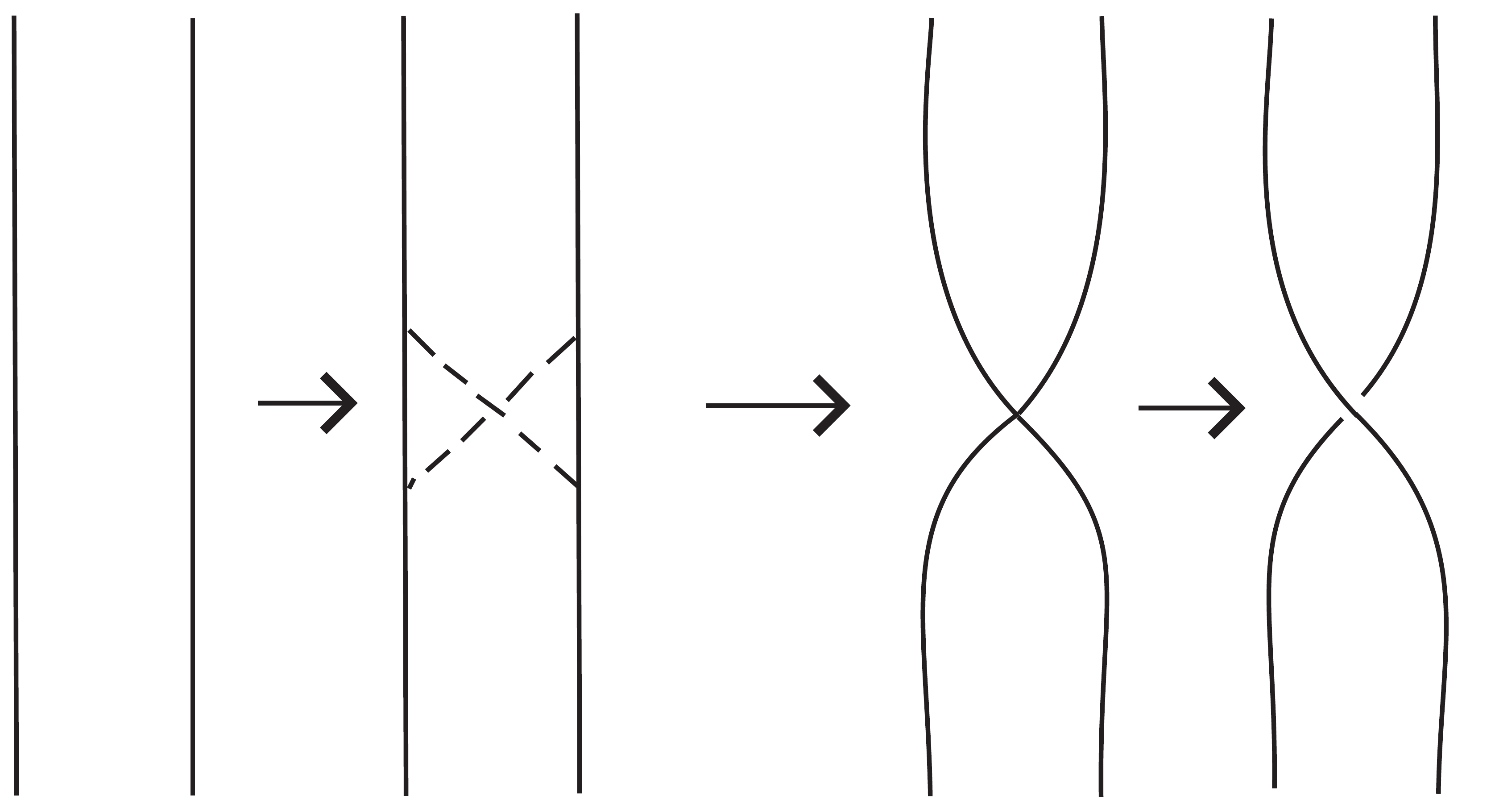}
\end{center}

\begin{center}
{\Huge $\Downarrow$}
\end{center}
\begin{center}
\noindent\llap{\raise 30pt\hbox{(2)}\qquad} \includegraphics[width=10cm]{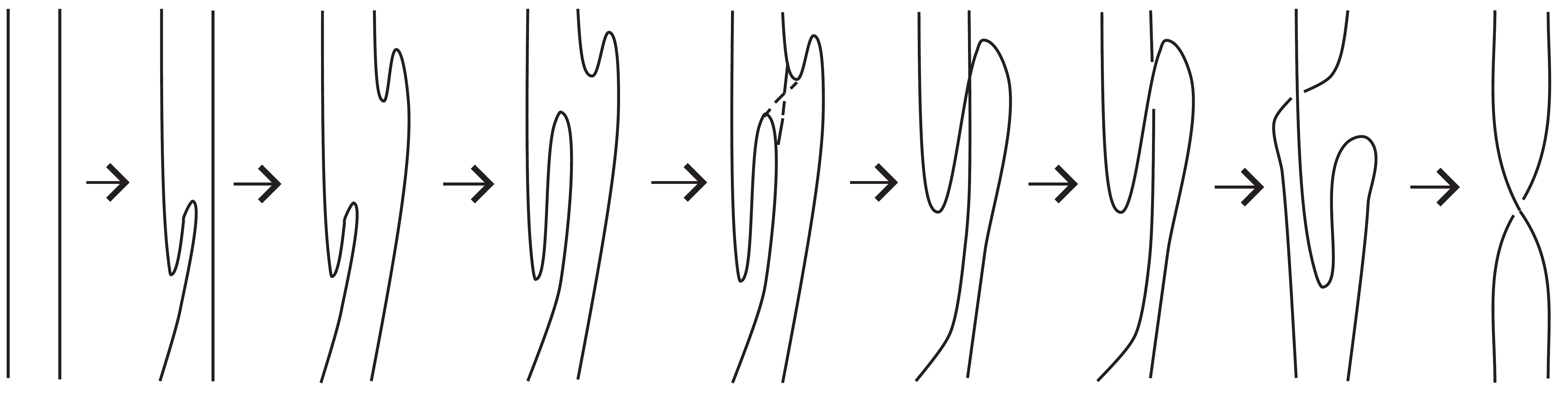}
\end{center}
\vspace{3mm}

\begin{center}
\noindent\llap{\raise 30pt\hbox{(3)}\qquad}\includegraphics[width=10cm]{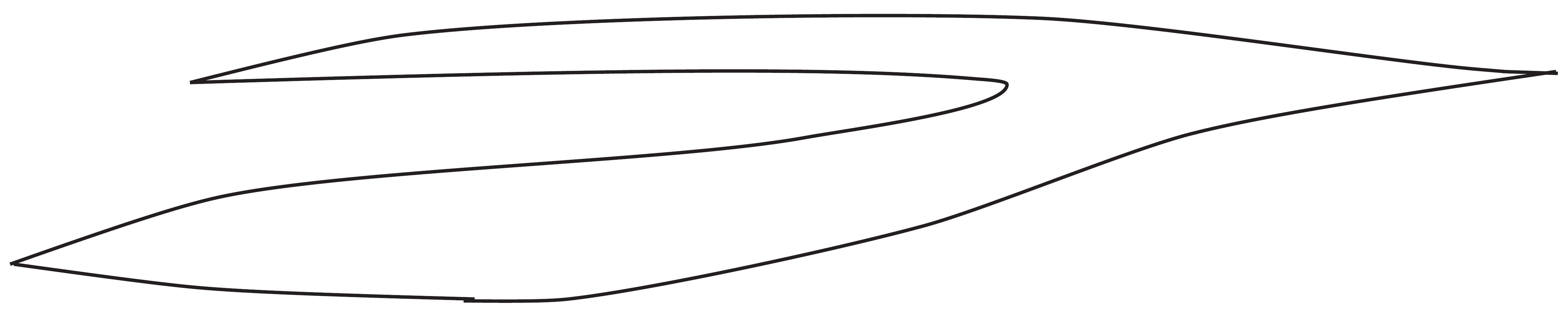}
\end{center}
\caption{Deformation of a branch point} 
\label{fig:q}
\vspace{10pt}
\end{figure}

\begin{figure}[ht!]
\begin{center}
\labellist\small
\pinlabel {$b=$Braid($F$)} [l] at 688 830
\pinlabel {$2b{+}2$} [b] <-5pt,0pt> at 223 479
\pinlabel {$2b{+}2$} [b] <5pt,0pt> at 373 465
\pinlabel {$2b{+}2$} [t] <-5pt,0pt> at 158 205
\pinlabel {$2b{+}2$} [t] <5pt,0pt> at 291 205 
\pinlabel {$2b$} at 0 389
\pinlabel {$2b{-}2$} [b] at -330 843
\pinlabel {$2$} at -280 490
\pinlabel {$0$} at -360 535
\endlabellist
\includegraphics[width=7cm]{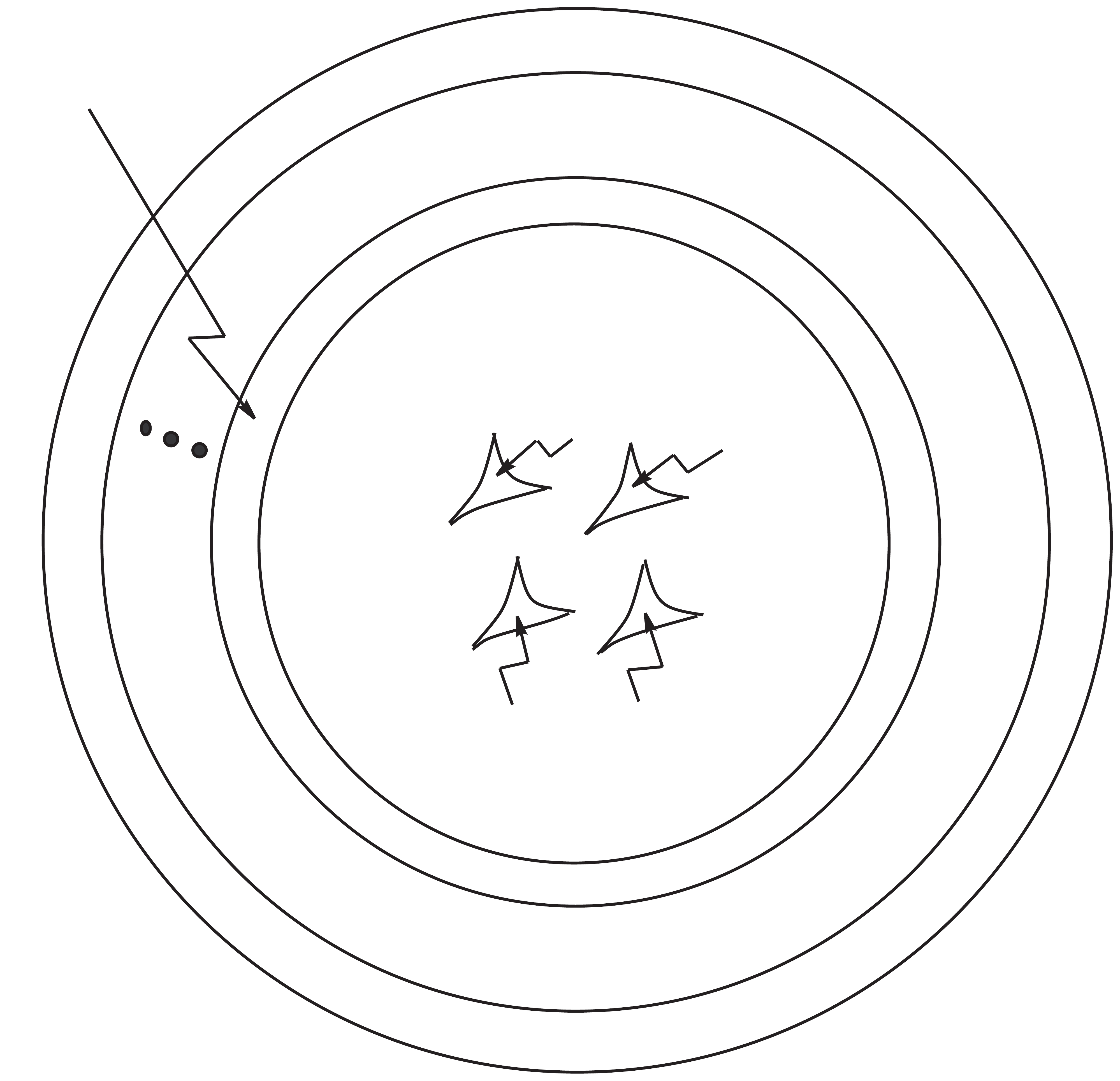}
\end{center}
\caption{Image of $S(\pi|_{\tilde{F}})$ by $\pi$ and local widths} 
\label{fig:r}
\vspace{10pt}
\end{figure}

\begin{proof}
Let $S^2 \subset \R^4$ be the 
standard $2$--sphere, ie, $S^2=\{(x_{1},x_{2},x_{3},x_{4})
\in \R^4 \ | \ x_{1}=0,\, x^2_{2}+x^2_{3}+x^2_{4}=1\}$,
and $D^2 \times S^2$ be its tubular neighborhood.
We may assume that $F \subset D^2 \times S^2$ and the restriction 
$pr_{2}|_{F}$ of $pr_{2} \co D^2 \times S^2 \to S^2$
is a simple branched covering of degree equal to Braid($F$).
We may further assume that the critical values of $pr_{2}|_{F}$
all lie near $(0,1,0,0) \in S^2$
and that outside of the pre-image of a neighborbhood of $(0,1,0,0)$,
$F$ is almost parallel to $S^2$.
Let us define the orthogonal projection $\pi \co \R^4 \to \R^2$
by $\pi(x_{1},x_{2},x_{3},x_{4})=(x_{3},x_{4})$.
Then, we may assume that the image of the singular points of
$\pi|_{F}$ is as depicted in \fullref{fig:p}.

Let $y \in S^2$ be a branch point of $g = pr_{2}|_{F}$ and let $x \in F$
be the branch point such that $y =g(x)$.
Furthermore, let $B \cong I \times J$ be a small neighborhood of $y$ in $S^2$
, where $I=J=[-1,1]$ and $y$ corresponds to (0,0),
and let $\tilde{B}$ be the component of $g^{-1}(B)$ which contains $x$.
Set $J_{t} = \{t\} \times J \subset I \times J$ for $t \in I$.
Then $(g|_{\tilde{B}})^{-1}(J_{t}) \subset pr^{-1}_{2}(J_{t})
\cong D^2 \times J$ 
can be regarded as a $2$--string braid for $t \neq 0$. See \fullref{fig:q} (1).

Then we deform $F$ (or more precisely, we deform $\tilde{B}$)
by an isotopy in $\R^4$ so that this sequence of $2$--string braids
is deformed as in \fullref{fig:q} (2).
Note that then $\pi$ is generic on $\tilde{B}$ and the image of the
singular points in $B \cong I \times J$ is as depicted in \fullref{fig:q} (3).
Three cusps are created, while the branch point in question
is eliminated.

We perform the above described deformation for each branch point of $g$.
Then we get a surface $\tilde{F}$ isotopic to $F$ such that $\pi$ is
generic with respect to $\tilde{F}$ and that the singular values of
$\pi|_{\tilde{F}}$ and the local widths are as depicted in \fullref{fig:r},
where $b=$Braid($F$).
Therefore, we have $w(F) \leq 2(b+1)$.
This completes the proof.
\end{proof}

\begin{figure}[ht!]
\labellist\small
\pinlabel {$b=$Braid($F$)} at 744 898
\pinlabel {$2b$} at 129 586
\pinlabel {$2b$} [l] at 361 451
\pinlabel {$2b$} [l] at 395 260
\pinlabel {$2b$} [r] at 80 473
\pinlabel* {$2b$} [t] at 111 221
\pinlabel {$2b{-}2$} [t] at 133 -200
\pinlabel {$2b{-}2$} [b] at 455 959
\pinlabel {$2$} [l] <0pt,2pt> at 795 406
\pinlabel {$2$} <1pt,5pt> at -282 489
\pinlabel {$0$} <0pt,5pt> at -360 508
\endlabellist
\begin{center}
\begin{minipage}{8.5cm}
\includegraphics[width=8.5cm]{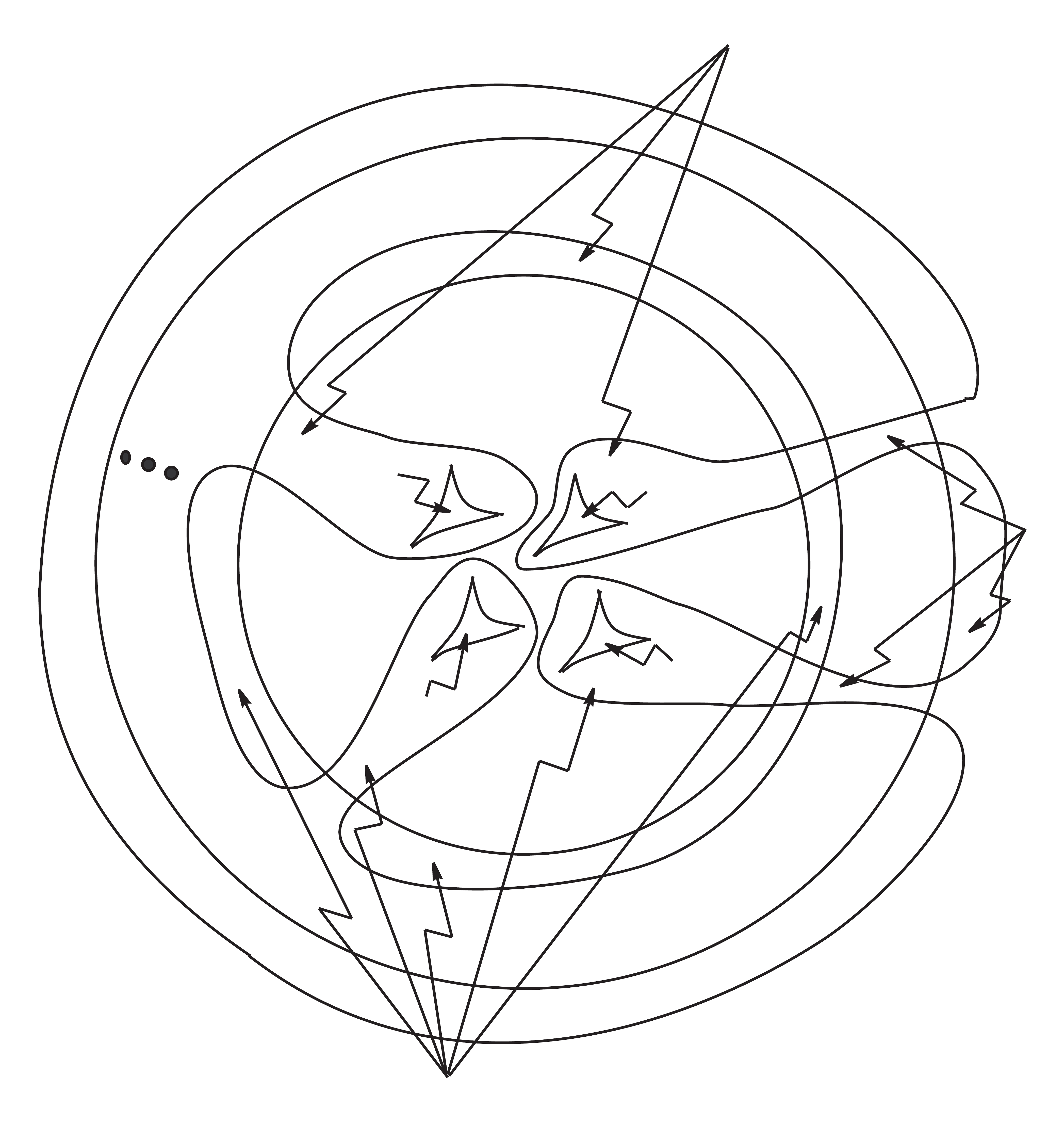}
\end{minipage}
\end{center}
\caption{A possible deformation of neighborhoods of the
corresponding branch points in $F$} 
\label{fig:r2}
\end{figure}

Let us consider a branch point of a surface braid as above.
Since it is simple, there may be a ``sheet" of $F$ over that
point which does not intersect a neighborhood of the
corresponding branch point in $F$.
If the sheet can be deformed as depicted in \fullref{fig:r2},
then the width decreases by $2$.
Therefore, the following conjecture seems to be plausible.
 
\begin{conjecture}
Let $F \subset \R^4$ be an orientable surface knot.
Then we have
\begin{eqnarray}
 w(F) \leq 2{\rm Braid}(F).   \nonumber
\end{eqnarray}
\end{conjecture}

By the following proposition, the difference between
the width and (twice) the braid index can be arbitrarily large.

\begin{proposition}\label{prop:t9}
For every $n \geq 3$, there exists a surface knot $F$ in $\R^4$ 
with {\rm Braid($F$)} $=n$ and $w(F)=4$.
\end{proposition} 

For the proof, we need the following.

\begin{lemma}\label{lem:t2}
For surface knots $F_{1}$ and $F_{2}$ in $\R^4$, we always have 
\begin{eqnarray}
 w(F_{1} \sharp F_{2}) \leq {\rm max}\{w(F_{1}),
w(F_{2})\}.   \nonumber
\end{eqnarray}
\end{lemma}

\begin{proof}
We may assume that there exists an orthogonal projection
$\pi \co \R^4 \to \R^2$ which is generic with respect to both
$F_{1}$ and $F_{2}$ such that $w(F_{1}, \pi) = w(F_{1})$ and
$w(F_{2}, \pi) = w(F_{2})$.
We may further assume that $\pi(F_{1}) \cap \pi(F_{2}) = \emptyset$.
Let us consider fold points of $\pi|_{F_{1}}$ and $\pi|_{F_{2}}$
whose images by $\pi$ lie in the outermost boundaries of
$\pi(F_{1})$ and $\pi(F_{2})$ respectively.
If we perform the connected sum operation using small disk neighborhoods
of these fold points and by connecting $F_{1}$ and $F_{2}$ by an
appropriate cylinder (see the proof of \fullref{prop:t0}),
then $\pi$ is generic with respect to $F_{1} \sharp F_{2}$ and
$w(F_{1} \sharp F_{2}, \pi) = {\rm max}\{w(F_{1}, \pi), w(F_{2}, \pi)\}$.
Thus the conclusion follows.
This completes the proof. 
\end{proof} 

\begin{remark}
The referee kindly pointed out that there is an example
for which the equality does not
hold in \fullref{lem:t2} as follows.
By Viro \cite{Vi}, it is known that there exists a ribbon $2$--sphere knot
$F \subset \R^4$,
which is not strongly trivial,
such that $F \sharp P_{+}$ is isotopic to $P_{+}$, where
$P_{+}$ is the trivial
projective plane with normal Euler number $2$. 
Let $K \subset \R^4$ be a Klein bottle knot, which is strongly trivial,
such that $P_{+} \sharp P_{-}$ is isotopic to $K$, where
$P_{-}$ is the trivial projective plane with normal Euler number $-2$.
Then $F \sharp K$ is isotopic to $K$.
Since $w(F)=4$, $w(K)=2$ and $w(F \sharp K)=2$,
the equality does not hold in \fullref{lem:t2}.
However, we do not know such an example if both $F_{1}$ and $F_{2}$
are orientable. 
The author would like to thank the referee for pointing out this example.
\end{remark}

\begin{proof}[Proof of \fullref{prop:t9}]
Let $F_{1}$ be the spun ($2,p$)--torus knot, where $p$ is an odd integer
with $p \geq 3$.
Furthermore, let $F$ be the connected sum of $n-2$ copies of $F_{1}$.
Then by Tanaka \cite{Ta}, we have Braid($F$) $=n$.
On the other hand, since the ($2,p$)--torus knot is a $2$--bridge knot,
we have $w(F_{1})=4$ by \fullref{prop:t1}.
Then by \fullref{lem:t2}, we have $w(F) \leq 4$.
Since Braid($F$) $=n \geq 3$, $F$ is not strongly trivial,
and hence $w(F) > 2$ by \fullref{thm:t0}.
Therefore, we have $w(F)=4$.
This completes the proof. 
\end{proof}

\section{Total widths of surface knots}\label{sec-5}

In this section, we give several characterization
theorems of surface knots with small total widths. 

The following is an immediate consequence of \fullref{thm:t0}.

\begin{theorem}\label{thm:t1}
Let $F \subset \R^4$ be a surface knot.
Then $tw(F)=2$ if and only if it is strongly trivial.
\end{theorem}

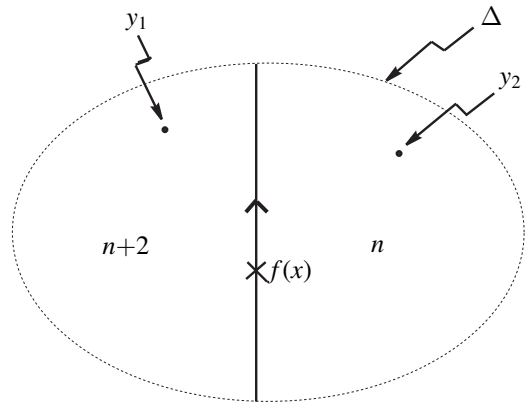
\begin{figure}[ht!]
\begin{center}
\labellist\small
\pinlabel {$n$} at 390 358
\pinlabel {$y_{2}$} [bl] at 535 547
\pinlabel {$\Delta$} [bl] at 511 628
\pinlabel {$y_{1}$} [b] at 90 619
\pinlabel {$n{+}2$} [r] at 115 358
\hair 1.5pt
\pinlabel {$f(x)$} [l] at 252 326
\endlabellist
\includegraphics[width=7cm]{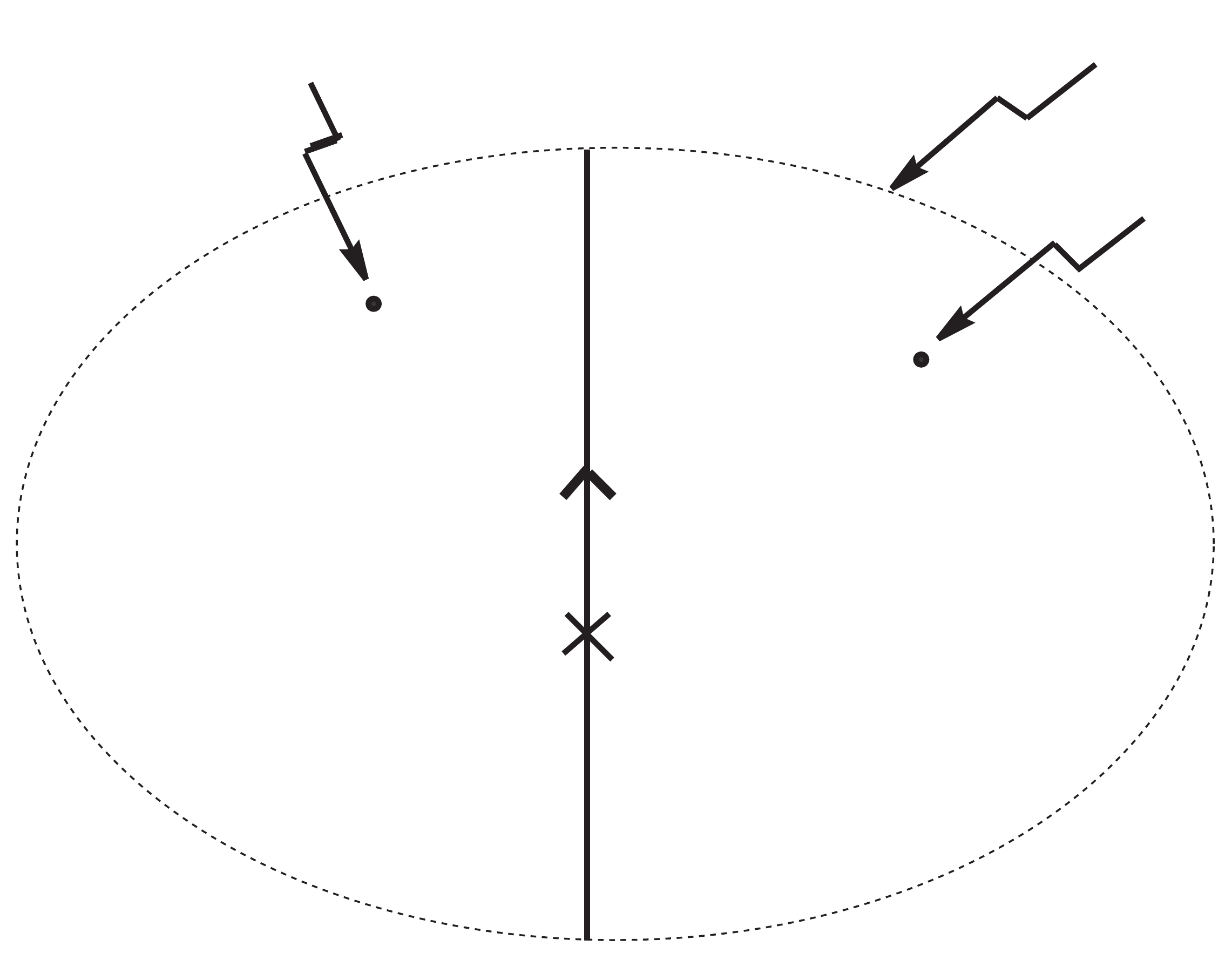}
\end{center}
\caption{Local orientation} 
\label{fig:s}
\vspace{10pt}
\end{figure}

Let $f \co F \to \R^2$ be a $C^{\infty}$ stable mapping of a
closed surface into the plane.
For a point $x \in S(f) \setminus S^2_{1}(f)$,
we give a local orientation of $S(f)$ at $x$ as follows.
For a sufficiently small disk neighborhood $\Delta$ of $f(x)$ in $\R^2$,
$\Delta \cap f(S(f))$ is an arc and $\Delta \setminus f(S(f))$
consists of two regions.
Let us take points, say $y_{1}$ and $y_{2}$, from each of the two
regions.
We may assume that the number of elements in the inverse image
$f^{-1}(y_{1})$ is greater than that of $f^{-1}(y_{2})$.
Then we orient $\Delta \cap f(S(f))$ so that the
left hand side region corresponds to $y_{1}$.
Finally we give a local orientation of $S(f)$ at $x$ so that
$f|_{S(f)}$ preserves the orientation around $x$.
See \fullref{fig:s}.

It is easy to see that the above local orientations vary
continuously and that they define a globally well-defined
orientation on $S(f)$.

On the other hand, by considering the ``line" $df_{x}(T_{x}S(f))$
for each $x \in S(f) \setminus S^2_{1}(f)$,
we obtain a smooth mapping $S(f) \setminus S^2_{1}(f) \to \RP^1$.
It is not difficult to see that this mapping extends to a 
smooth mapping $\tau_{f} \co S(f) \to \RP^1$.
We orient $\RP^1$ so that the lines rotating in the counter-clockwise
direction correspond to the positive direction of $\RP^1$.

Then we define rot($f$) to be the mapping degree of 
$\tau_{f} \co S(f) \to \RP^1$.

Then the following lemma is proved in Levine \cite{Le}.

\begin{lemma}\label{lem:le}
The Euler characteristic $\chi(F)$ of $F$ coincides with
{\rm rot(}$f${\rm)}.
\end{lemma}

Using \fullref{lem:le}, we prove the following.

\begin{theorem}\label{thm:t2}
Let $F \subset \R^4$ be a surface knot which is diffeomorphic 
to the $2$--sphere $S^2$.
Then
$tw(F) \leq 6$ if and only if it is strongly trivial.
\end{theorem}

\begin{figure}[ht!]
\cl{\includegraphics[width=9cm]{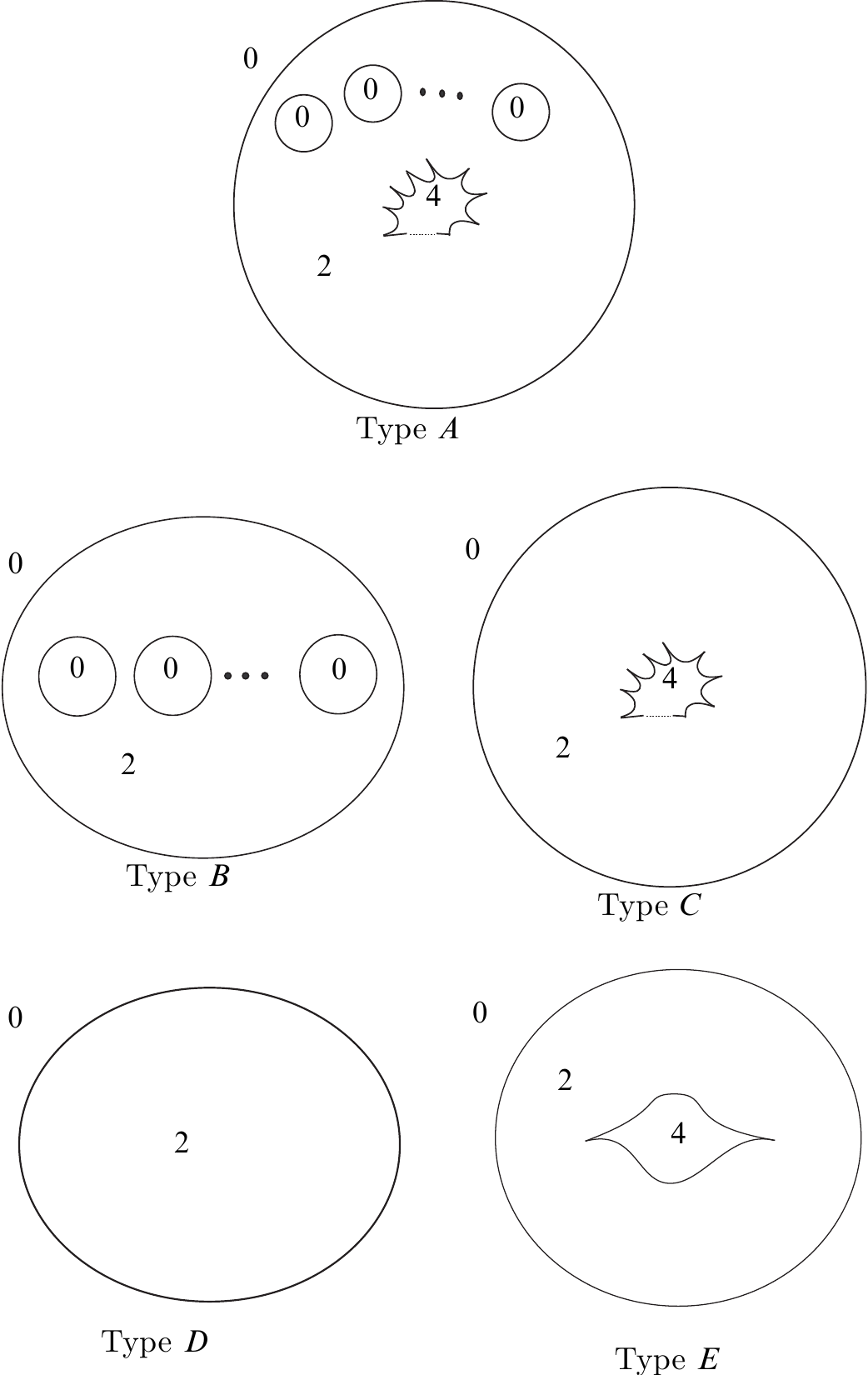}}
\caption{Possible images of the singular set with no fold crossing for
$F \cong  S^2$} 
\label{fig:t}
\end{figure}

\begin{figure}[ht!]
\cl{\includegraphics[width=9cm]{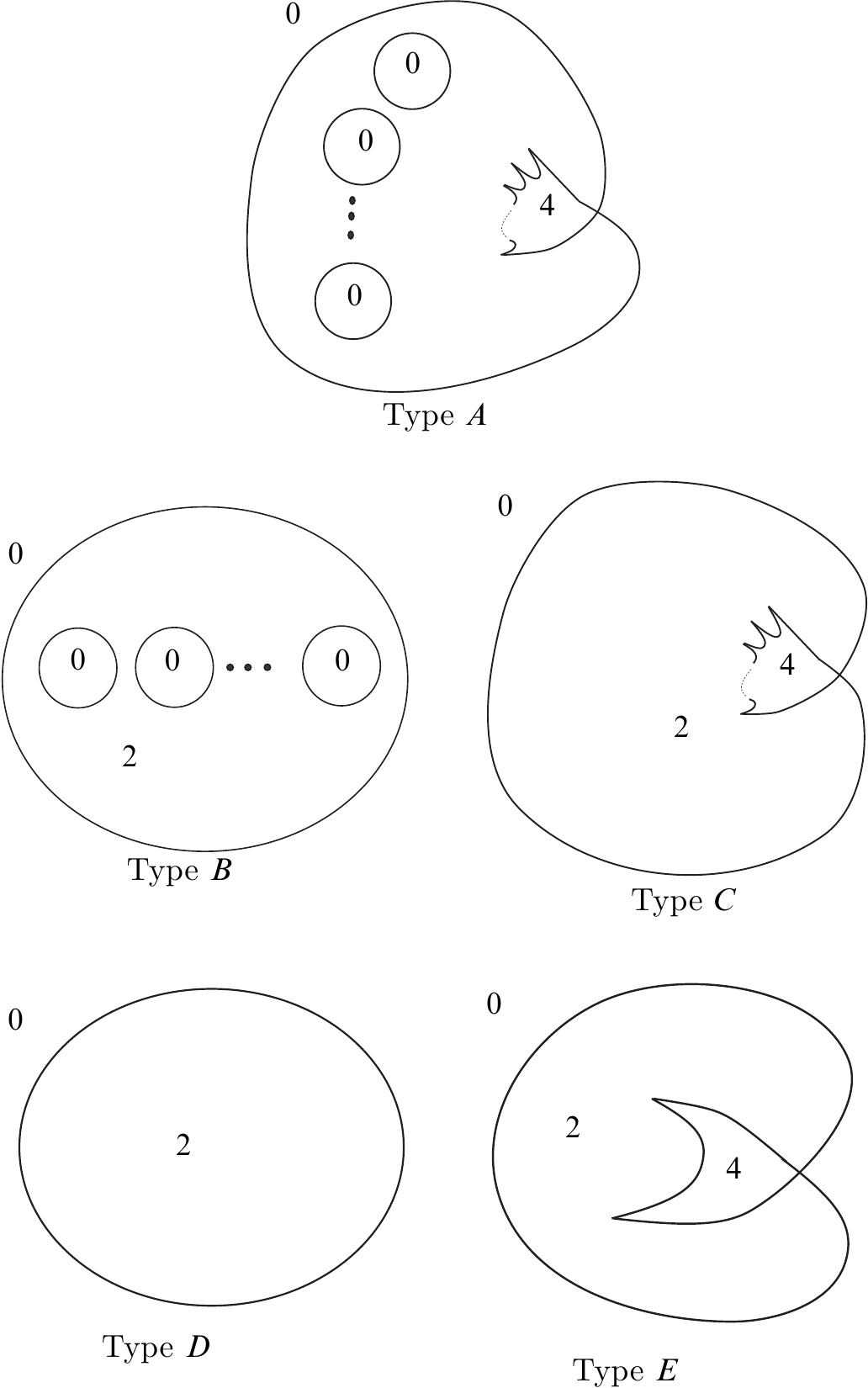}}
\caption{Possible images of the singular set with one fold crossing for
$F \cong  S^2$,
part~$1$} 
\label{fig:u}
\vspace{10pt}
\end{figure}

\begin{figure}[ht!]
\cl{\includegraphics[width=9cm]{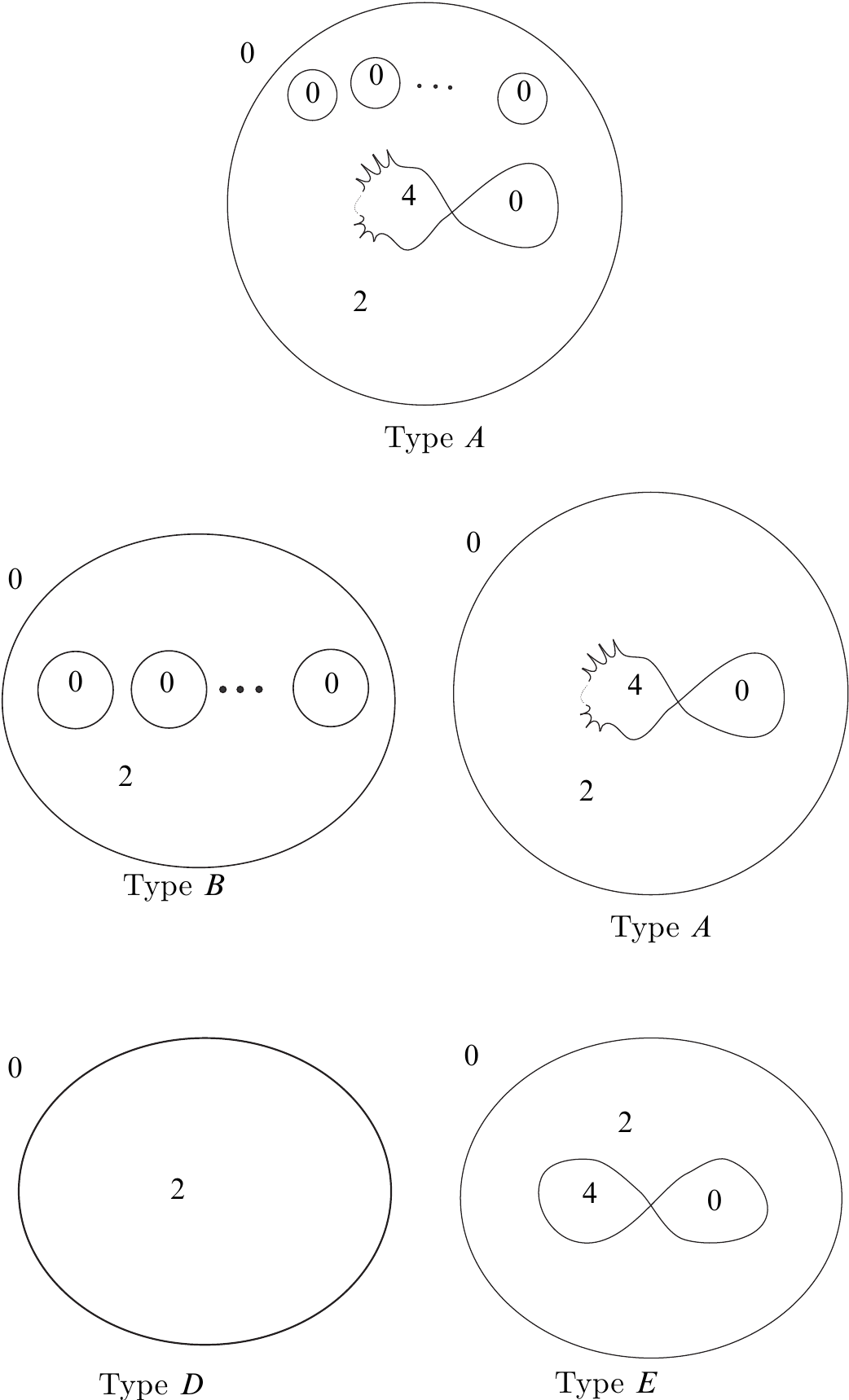}}
\caption{Possible images of the singular set with one fold crossing for
$F \cong  S^2$,
part~$2$} 
\label{fig:v}
\vspace{10pt}
\end{figure}

\begin{figure}[ht!]
\cl{\includegraphics[width=9cm]{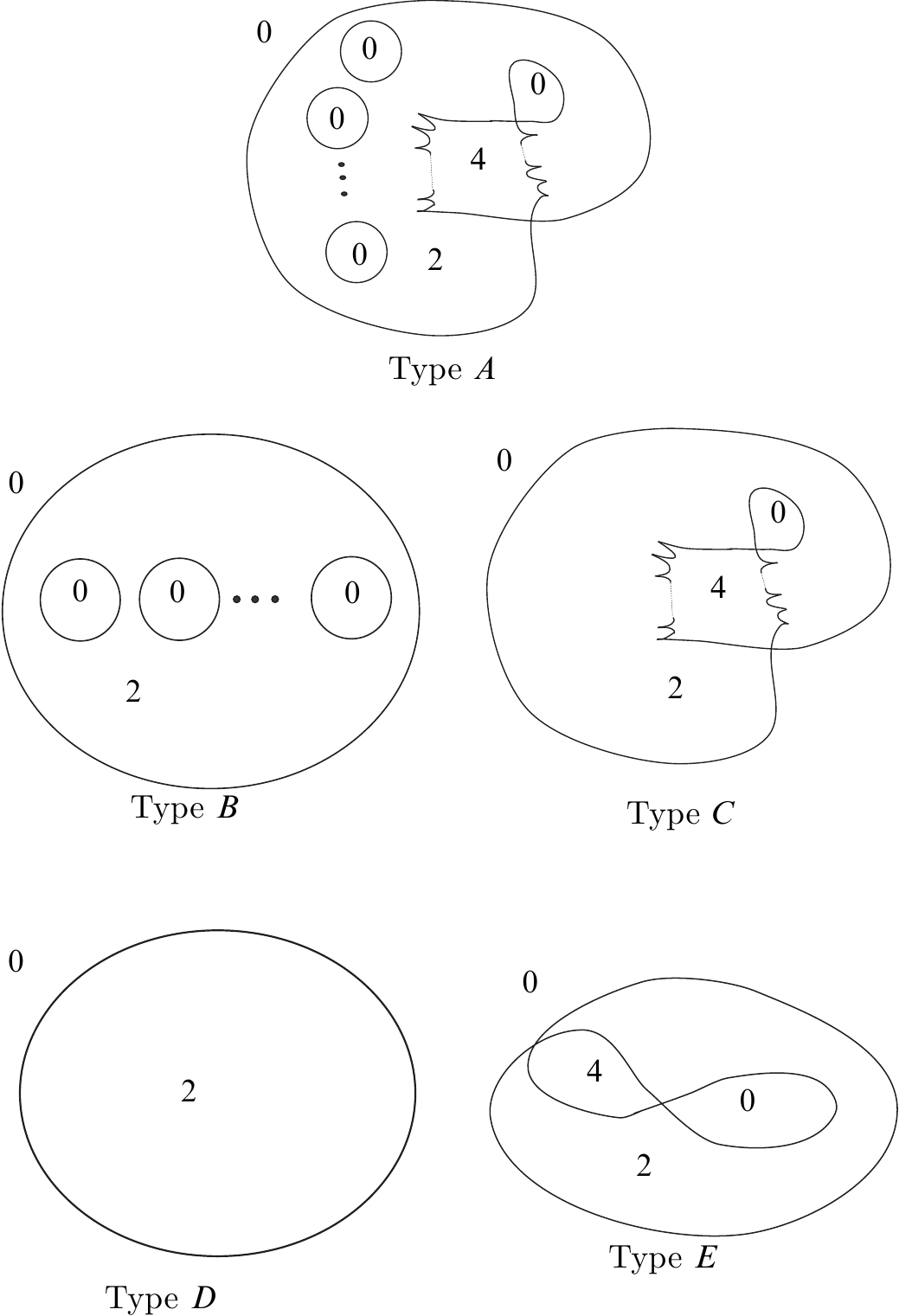}}
\caption{Possible images of the singular set with two fold crossings for 
$F \cong  S^2$, part~$1$} 
\label{fig:w}
\end{figure}

\begin{figure}[ht!]
\cl{\includegraphics[width=5.5cm]{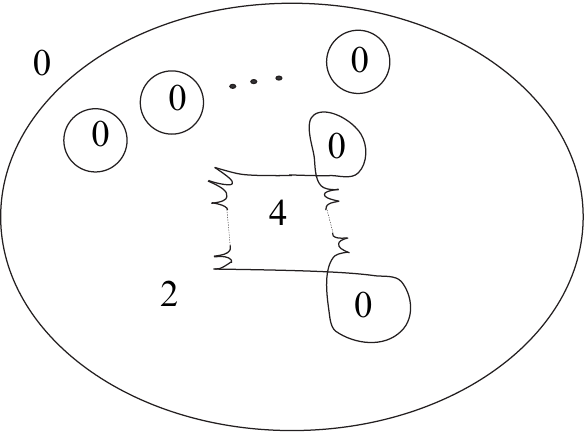}}
\caption{Possible image of the singular set with two fold crossings
for $F \cong  S^2$, part~$2$} 
\label{fig:x}
\vspace{20pt}
\end{figure}

\begin{figure}[ht!]
\cl{\includegraphics[width=11cm]{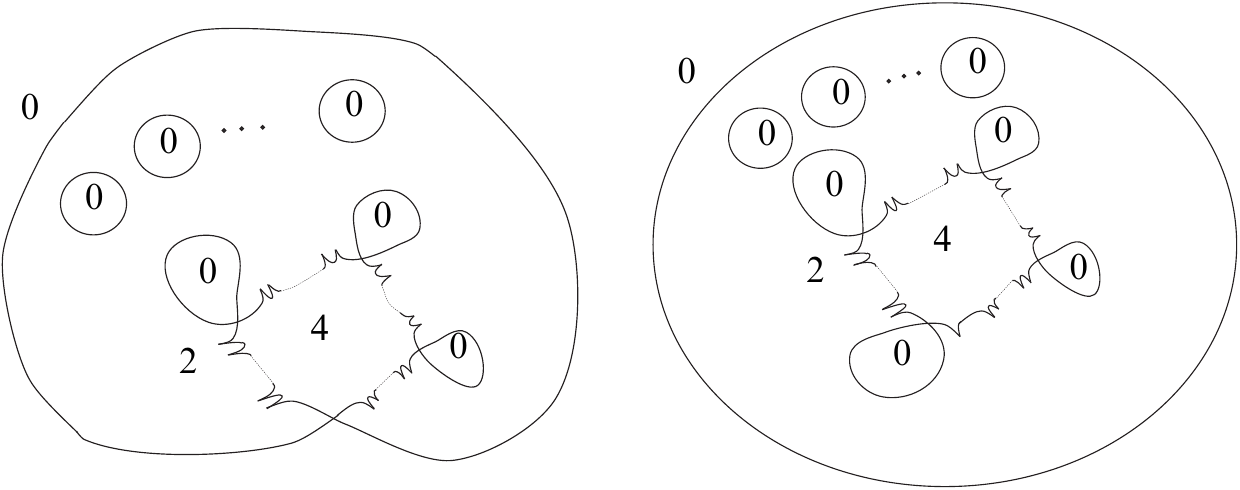}}
\caption{Possible images of the singular set with three or more fold 
crossings for $F \cong  S^2$} 
\label{fig:y}
\vspace{20pt}
\end{figure}

\begin{proof}
If $tw(F)=2$, then by \fullref{thm:t1}, $F$ is strongly trivial.
Furthermore, there does not exist a surface knot 
$F$ with $tw(F)=4$, since $F$ is connected.
Therefore, we may assume $tw(F)=6$ and there exists an orthogonal projection 
$\pi \co \R^4 \to \R^2$ which is generic with respect to $F$
such that $tw(F, \pi)=tw(F)$.

If $\pi(S(\pi|_{F}))$ has no fold crossings, then it is of the form
``Type~$A$" as depicted in \fullref{fig:t} up to isotopy of $\R^2$.
Then, by \fullref{lem:t0}, $F$ is the connected sum of surface knots
$F_{1}$ and $F_{2}$ such that $\pi(S(\pi|_{F_{1}}))$ 
(or $\pi(S(\pi|_{F_{2}}))$) is of the form ``Type~$B$" (resp.\ ``Type~$C$")
as depicted in \fullref{fig:t}
up to isotopy of $\R^2$.
Since $F$ is diffeomorphic to the $2$--sphere, so are $F_{1}$ and $F_{2}$.
Then, by \fullref{lem:le}, Type~$B$ and Type~$C$ must correspond to
Type~$D$ and Type~$E$ of \fullref{fig:t} respectively.
By \fullref{lem:crs}, we see that $F_{1}$ is strongly trivial.
Furthermore, there exists an orthogonal projection ${\displaystyle \pi^4_{1}}
\co \R^4 \to \R^1$ which is generic with respect to $F_{2}$ such that
${\displaystyle \pi^4_{1}}|_{F_{2}}$ has exactly two critical points.
In fact, such a projection can be obtained by composing $\pi \co \R^4 \to \R^2$
and a suitable projection $\R^2 \to \R^1$ (for example, see Fukuda \cite{F}).
Thus, $F_{2}$ is also strongly trivial, and hence so is $F=F_{1} \sharp F_{2}$.

If $\pi(S(\pi|_{F}))$ has one fold crossing, then it is of the form
``Type~$A$" as depicted in \fullref{fig:u} or in 
\fullref{fig:v} up to isotopy
of $\R^2$. Then, by \fullref{lem:t0}, $F$ is the connected
sum of surface knots
$F_{1}$ and $F_{2}$ such that $\pi(S(\pi|_{F_{1}}))$
(or $\pi(S(\pi|_{F_{2}}))$) is of the form ``Type~$B$" (resp.\ ``Type~$C$")
as depicted in \fullref{fig:u}
or in \fullref{fig:v} up to isotopy of $\R^2$.
Since $F$ is diffeomorphic to the $2$--sphere, so are $F_{1}$ and $F_{2}$.
Then, by \fullref{lem:le}, Type~$B$ and Type~$C$ must correspond to
Type~$D$ and Type~$E$ of \fullref{fig:u} or \fullref{fig:v} respectively.
By \fullref{lem:crs}, $F_{1}$ is strongly trivial.
Furthermore, there exists an orthogonal projection
${\displaystyle \pi^4_{1}} \co \R^4 \to \R^1$ which is generic with respect to
$F_{2}$ such that ${\displaystyle \pi^4_{1}}|_{F_{2}}$ has 
exactly four critical points.
Therefore, $F_{2}$ is strongly trivial by Scharlemann \cite{Sc}.
(In fact, ``Type~$E$" of \fullref{fig:v} does not occur by
Akhmet'ev \cite[23.~Corollary]{Ak}.)
Thus, $F=F_{1} \sharp F_{2}$ is strongly trivial.

If $\pi(S(\pi|_{F}))$ has two fold crossings,
then it is of the form ``Type~$A$" as depicted in \fullref{fig:w}
or as depicted in \fullref{fig:x}.
In the former case, we see that $F$ is strongly trivial as before
(see \fullref{fig:w}).
In the latter case, we see that $\chi(F) < \chi(S^2)$ by \fullref{lem:le},
which is a contradiction.
Thus, this case does not occur.

If $\pi(S(\pi|_{F}))$ has three or more fold crossings,
then it is of the form as depicted in \fullref{fig:y}.
Then, we see that $\chi(F) < \chi(S^2)$ by \fullref{lem:le},
so that this case does not occur.

Hence $F$ is always strongly trivial.
This completes the proof.  
\end{proof}

\begin{corollary}
Let $F \subset \R^4$ be an $n$--twist spun 2--bridge knot with $n \neq \pm 1$.
Then we have $tw(F)=8$.
\end{corollary}

\begin{proof}
Since $F$ is not strongly trivial, by \fullref{thm:t2}
we have $tw(F) \geq 8$.
On the other hand, since $F$ has planar projection as in \fullref{fig:o},
we have $tw(F) \leq 8$.
This completes the proof.
\end{proof}

Similarly, for surface knots diffeomorphic to the projective plane,
we have the following characterization.

\begin{theorem}
Let $F$ be a surface knot which is diffeomorphic to the
projective plane $\RP^2$. Then
$tw(F) \leq 6$ if and only if it is trivial.
\end{theorem}

\begin{figure}[ht!]
\begin{center}
\includegraphics[width=9cm]{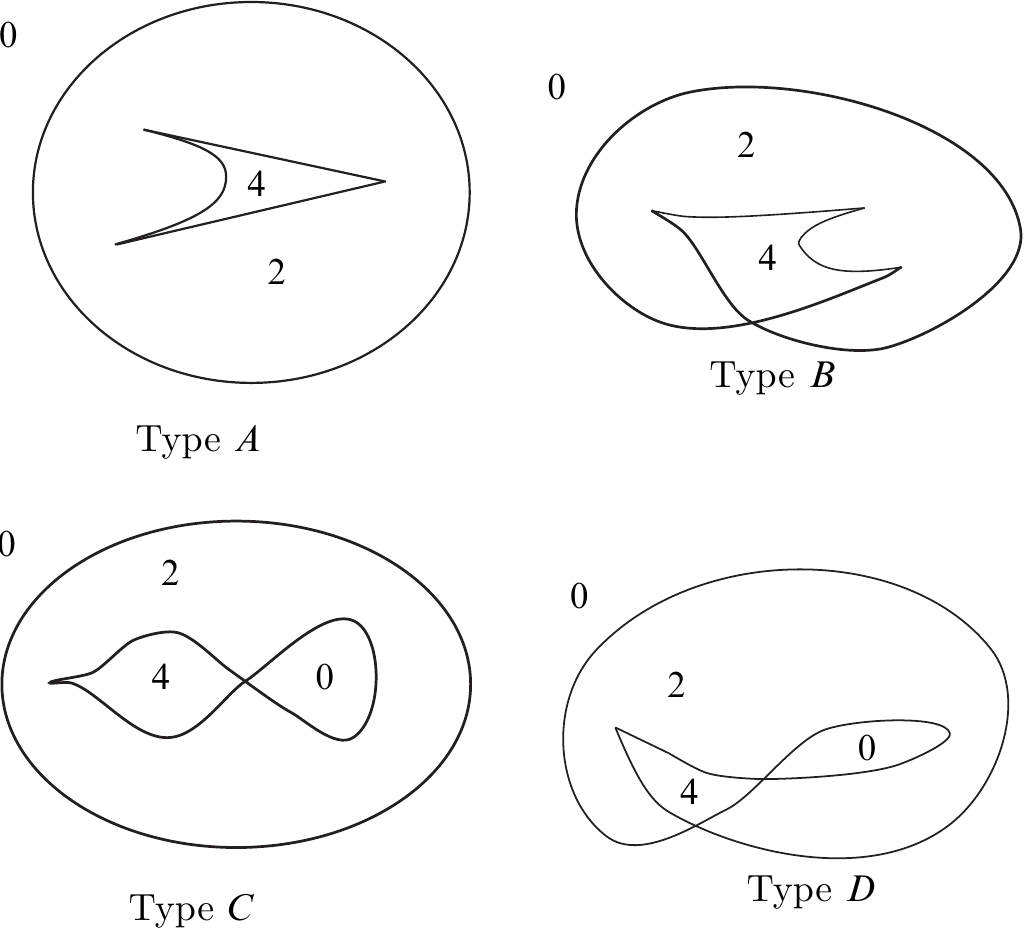}
\end{center}
\caption{Possible images of the singular set for $F \cong \RP^2$ with 
$tw(F)=6$} 
\label{fig:z}
\end{figure}

\begin{proof}
If $tw(F)=2$, then by \fullref{thm:t1}, $F$ is strongly trivial.
Furthermore, there does not exist a surface knot $F$ with $tw(F)=4$,
since $F$ is connected.
Therefore, we may assume $tw(F)=6$ and there exists an orthogonal projection 
$\pi \co \R^4 \to \R^2$ which is generic with respect to $F$
such that $tw(F, \pi)=tw(F)$.

We use the argument of the proof of \fullref{thm:t2}. 
If $\pi(S(\pi|_{F}))$ has no fold crossings, then it is of the form
``Type~$A$" as depicted in \fullref{fig:t}.
Since $F$ is diffeomorphic to the projective plane,
by \fullref{lem:t0} we see that $F=F_{1} \sharp F_{2}$, where 
$\pi(S(\pi|_{F_{1}}))$ is of the form ``Type~$D$" as depicted in
\fullref{fig:t} and $\pi(S(\pi|_{F_{2}}))$ is of the form ``Type~$A$"
as depicted in \fullref{fig:z}.
By \fullref{lem:crs}, $F_{1}$ is strongly trivial.
Since there exists an orthogonal projection 
${\displaystyle \pi^4_{1}} \co \R^4 \to \R^1$
which is generic with respect to $F_{2}$ such that
${\displaystyle \pi^4_{1}}|_{F_{2}}$ has exactly three critical points,
we see that $F_{2}$ is trivial by Bleiler and Scharlemann \cite{bs}.
Therefore, $F$ is trivial.

If $\pi(S(\pi|_{F}))$ has one fold crossing, then it is of the form
``Type~$A$" as depicted in \fullref{fig:u} or in 
\fullref{fig:v} up to isotopy of $\R^2$. 
Since $F$ is diffeomorphic to the projective plane,
by \fullref{lem:t0} we see that $F=F_{1} \sharp F_{2}$,
where $\pi(S(\pi|_{F_{1}}))$ is of the form ``Type~$D$" as depicted
in \fullref{fig:u} or \fullref{fig:v} and $\pi(S(\pi|_{F_{2}}))$
is of the form ``Type~$B$" or ``Type~$C$" as depicted in \fullref{fig:z}.
By \fullref{lem:crs}, $F_{1}$ is strongly trivial.
Since there exists an orthogonal projection 
${\displaystyle \pi^4_{1}} \co \R^4 \to \R^1$
which is generic with respect to $F_{2}$ such that
${\displaystyle \pi^4_{1}}|_{F_{2}}$ has exactly three critical points,
we see that $F_{2}$ is trivial by \cite{bs}.
Therefore, $F$ is trivial.

If $\pi(S(\pi|_{F}))$ has two fold crossings,
then it is of the form ``Type~$A$" as depicted in \fullref{fig:w}
or as depicted in \fullref{fig:x}.
In the former case, we see that $F$ is trivial as before
(see ``Type~$D$" in \fullref{fig:z}).
In the latter case, we see that $\chi(F) < \chi(\RP^2)$ by \fullref{lem:le},
which is a contradiction.
Thus, this case does not occur.

If $\pi(S(\pi|_{F}))$ has three or more fold crossings,
then it is of the form as depicted in \fullref{fig:y}.
Then, we see that $\chi(F) < \chi(\RP^2)$ by \fullref{lem:le},
so that this case does not occur.

Hence $F$ is always trivial.
This completes the proof.  
\end{proof}

\begin{figure}[ht!]
\begin{center}
\includegraphics[width=9cm]{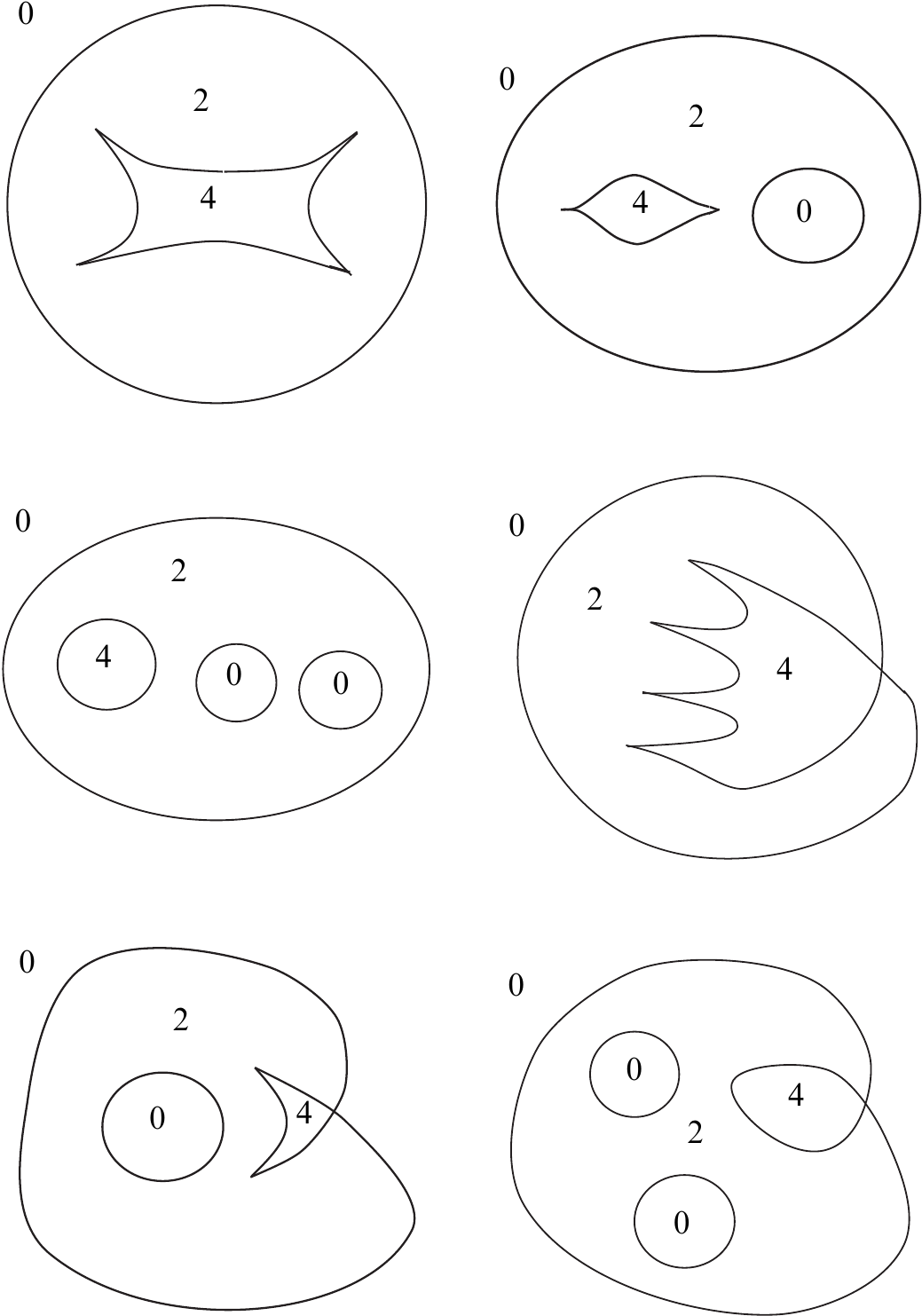}
\end{center}
\caption{Possible configurations for the Klein bottle case, part~$1$} 
\label{fig:ab}
\end{figure}

\begin{figure}[ht!]
\begin{center}
\includegraphics[width=9cm]{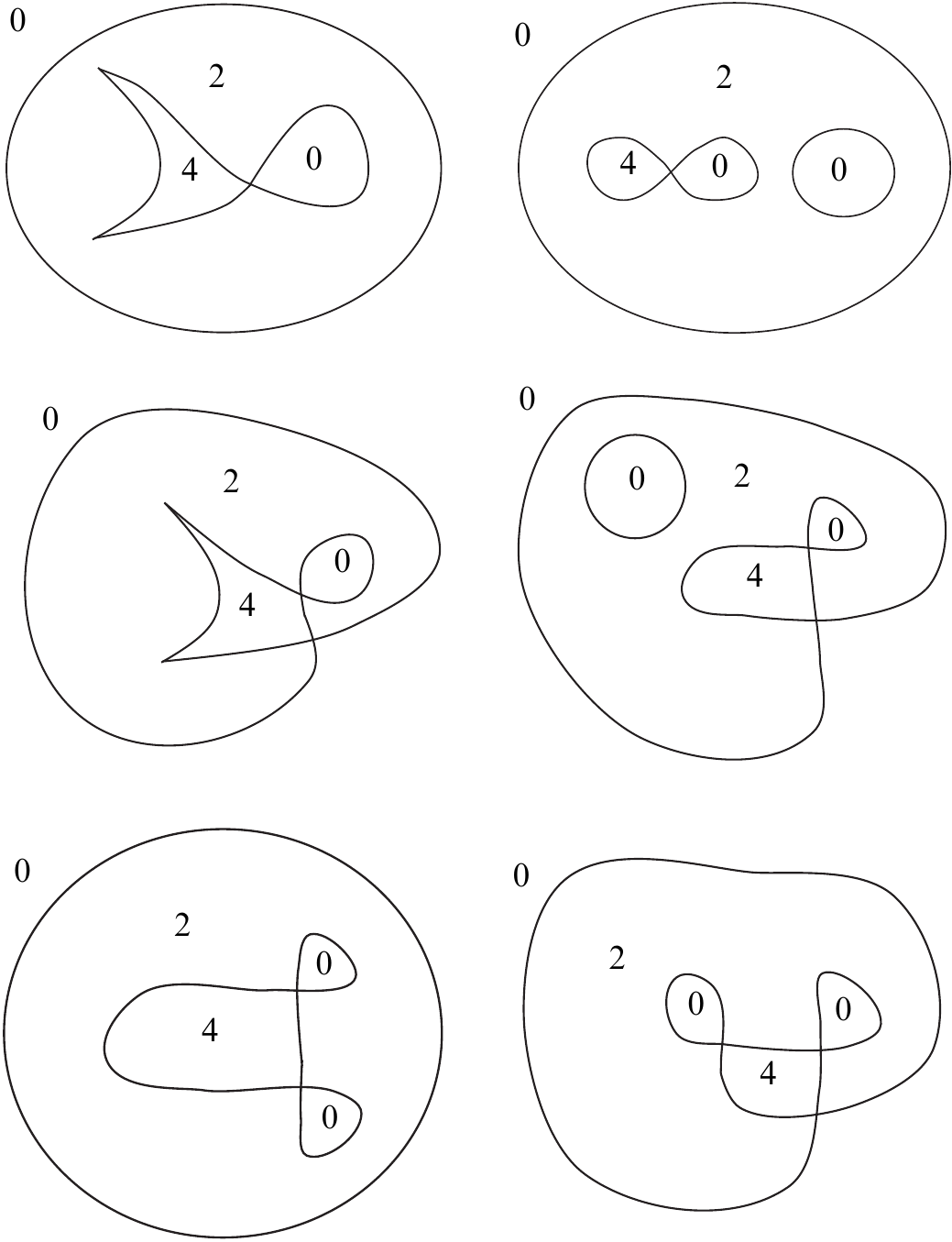}
\end{center}
\caption{Possible configurations for the Klein bottle case, part~$2$} 
\label{fig:ac}
\end{figure}

\begin{remark}
We do not know if a similar characterization theorem holds for surface
knots of higher genus.
For example, in Figures~\ref{fig:ab} and \ref{fig:ac} we have listed 
all the possible configurations of the planar image of the singular set
for knotted Klein bottles with total width smaller than or equal to six.
In general, we have many cusps and cannot apply \fullref{lem:crs}
directly.
Furthermore, we have no unknotting theorem as in \cite{bs,Sc}
for embedded Klein bottles as far as the author knows.
\end{remark}

\bibliographystyle{gtart}
\bibliography{link}

\end{document}